\documentclass[11pt]{article}
\usepackage{amsmath,amssymb,amsthm,textcomp}
\usepackage[usenames]{color}\usepackage{graphicx}
\usepackage{enumerate}
\definecolor{brown}{cmyk}{0, 0.72, 1, 0.45}
\definecolor{grey}{gray}{0.5}

\def\z{\zeta}

\allowdisplaybreaks
\newcommand{\rdown}[1]{{\left\lfloor #1\right \rfloor}}
\def\wel{1\wedge\l}

\def\HAM{{\sc extend-rotate}}
\def\cG{{\cal G}}
\newcommand{\card}[1]{\left|#1\right|}
\newcommand{\gap}[1]{\hspace{#1in}}
\def\sd{\d^*}
\def\ha{\hat{\alpha}}
\def\hb{\hat{\b}}
\def\hx{\hat{\xi}}
\def\hA{\hat{A}}
\def\hB{\hat{B}}
\def\hC{\hat{C}}
\def\hD{\hat{D}}

\def\bu{{\bf u}}

\def\lb{\bar{\l}}
\def\s{\sigma}
\def\si{\par\smallskip\noindent}

\def\ex{\mathbb{E\/}}

\def\la{\lambda}
\def\la{\lambda}
\def\a{\alpha}

\def\part{\partial}
\def\tag#1 {\eqno(#1)}

\newcommand{\brac}[1]{\left( #1 \right)}
\newcommand\bfrac[2]{\left(\frac{#1}{#2}\right)}
\newtheorem{Theorem}{Theorem}

\newtheorem{claim}{Claim}[Theorem]

\def\a{\alpha} \def\b{\beta} \def\d{\delta} \def\D{\Delta}
\def\e{\epsilon} \def\f{\phi}   
\def\G{\Gamma}  \def\k{\kappa} 
\def\z{\zeta} \def\th{\theta}   \def\l{\lambda}
\def\La{\Lambda} \def\m{\mu} \def\n{\nu} \def\p{\pi}
\def\r{\rho}  \def\s{\sigma} 
\def\t{\tau} \def\om{\omega}  

\def\bx{{\bf x}}
\def\by{{\bf y}}
\def\bz{{\bf z}}
\def\bd{{\bf d}}
\def\bv{{\bf v}}
\def\bw{{\bf w}}
\def\sx{{\bf x}^*}

\def\M{M^*}
\def\N{N^*}

\def\tiy{\tilde{y}}
\def\tiz{\tilde{z}}
\def\tim{\tilde{\m}}
\def\til{\tilde{\l}}
\def\tiT{\tilde{T}}
\def\tiu{\tilde{u}}
\def\tiv{\tilde{v}}

\def\hhu{\hat{u}}
\def\hhv{\hat{v}}

\def\cE{{\cal E}}
\def\cD{{\cal D}}
\def\e{\varepsilon}
\def\n{\nu}

\newtheorem{lemma}{Lemma}[section]
\newcommand{\proofstart}{{\bf Proof\hspace{2em}}}
\newcommand{\proofend}{\hspace*{\fill}\mbox{$\Box$}}
\newcommand{\set}[1]{\left\{#1\right\}}
\newcommand{\ignore}[1]{}

\def\hy{\hat{y}}
\def\hz{\hat{z}}
\def\hu{\hat{\bu}}

\newcommand{\beq}[1]{\begin{equation}\label{#1}}
\def\eeq{\end{equation}}
\def\hm{\hat{\m}}
\def\hl{\hat{\l}}
\def\t{\tau}\def\hT{\hat{T}}
\def\La{\Lambda}
\def\bb{{\bf b}}
\def\whp{\text{w.h.p.}}
\def\2G{{\sc 2greedy}}

\def\Pr{\mathbb{P}}
\def\E{\mathbb{E}}

\def\cP{{\cal P}}
\def\qs{\text{q.s.}}

\begin{document}

\title{On a Greedy 2-Matching Algorithm and Hamilton Cycles in Random Graphs with Minimum Degree at Least Three.}

\author{Alan Frieze\thanks{Research supported in part by
NSF Grant CCF1013110},\\Department of Mathematical Sciences,\\
Carnegie Mellon University,\\Pittsburgh PA15217.}

\maketitle
\begin{abstract}
We describe and analyse a simple greedy algorithm \2G\ that finds a good 2-matching $M$ in
the random graph $G=G_{n,cn}^{\d\geq 3}$ when $c\geq 15$. A 2-matching is a spanning subgraph
of maximum degree two and $G$ is drawn uniformly from graphs with vertex set $[n]$, $cn$ edges and
minimum degree at least three. By good we mean that $M$ has $O(\log n)$ components. We then use this
2-matching to build a Hamilton cycle in $O(n^{1.5+o(1)})$ time \whp.
\end{abstract}

\section{Introduction}
There have been many papers written on the existence of Hamilton cycles in random graphs.
Koml\'os and Szemer\'edi \cite{KoSz}, Bollob\'as \cite{B1}, Ajta, Koml\'os and Szemer\'edi \cite{AKS}
showed that the question is intimately related to the minimum degree.
Loosely speaking, if we are considering random graphs with
$n$ vertices and minimum degree at least two then we need $\Omega(n\log n)$
edges in order that they are likely to be Hamiltonian.

For sparse random graphs with $O(n)$
random edges, one needs to have minimum degree at least three. This
is to avoid having three vertices of degree two sharing a common neighbour.
There are several models of a random graph in which minimum degree three is satisfied:
Random regular graphs of degree at least
three, Robinson and Wormald \cite{RW1}, \cite{RW2} or
the random graph $G_{3-out}$, Bohman and Frieze \cite{BF}.
Bollob\'as, Cooper, Fenner and Frieze \cite{BCFF} considered
the classical random graph $G_{n,m}$ with conditioning on
the minimum degree $k$ i.e. each graph with vertex set $[n]$ and
$m$ edges and minimum degree at least $k$ is considered to be
equally likely.
Denote this model of a random graph by $G_{n,m}^{\d\geq k}$.
They showed that for every $k\geq 3$ there is a $c_k\leq (k+1)^3$
such that if $c\geq c_k$ then \whp\ $G_{n,cn}^{\d\geq k}$
has $(k-1)/2$ edge disjoint Hamilton cycles, where a perfect matching
constitutes half a Hamilton cycle in the case where $k$ is even.
It is reasonable to conjecture that $c_k=k/2$. The results of this paper
and a companion \cite{FP} reduce the known value
of $c_3$ from 64 to below 15. It can be argued that replacing one incorrect upper
bound by a smaller incorrect upper bound does not constitute significant progress. However,
the main contribution of this paper is to introduce a new
greedy algorithm for finding a large 2-matching in a random graph and to give a (partial) analysis
of its performance and of course to apply it to the Hamilton cycle problem.

One is interested in the time taken to construct a Hamilton cycle
in a random graph. Angluin and Valiant \cite{AV}
and Bollob\'as, Fenner and Frieze \cite{BFF} give polynomial time algorithms.
The algorithm in \cite{AV} is very fast, $O(n\log^2n)$
time, but requires $Kn\log n$ random edges for sufficiently large $K>0$. The
algorithm in \cite{BFF} is of order $n^{3+o(1)}$ but
works \whp\ at the exact threshold for connectivity. Frieze \cite{F1} gave an
$O(n^{3+o(1)})$ time algorithm for finding
large cycles in sparse random graphs and this could be adapted to find Hamilton
cycles in $G_{n,cn}^{\d\geq 3}$ in this time for sufficiently
large $c$. Another aim of \cite{FP} and this paper is reduce this running time.
The results of this paper and its companion \cite{FP}
will reduce this to $n^{1.5+o(1)}$ for sufficiently large $c$, and perhaps in a later
paper, we will further reduce the running time by
borrowing ideas from a linear expected time algorithm for matchings due to
Chebolu, Frieze and Melsted \cite{CFM}.

The idea of \cite{CFM} is to begin the process of constructing a
perfect matching by using the Karp-Sipser algorithm \cite{KS} to find
a good matching and then build this up to a perfect matching by alternating
paths. The natural extension of this idea is to
find a good 2-matching and then use extension-rotation arguments to
transform it into a Hamilton cycle.
A 2-matching $M$ of $G$ is a spanning subgraph of maximum degree 2.
Each component of $M$ is a cycle or a path (possibly an isolated
vertex) and we let $\k(M)$ denote the number of components of $M$.
The time taken to transform $M$ into a Hamilton cycle depends heavily on $\k(M)$.
The aim is to find a 2-matching $M$ for which $\k(M)$ is small.
The main result of this paper is the following:
\begin{Theorem}\label{th1}\
There is an absolute constant $c_0>0$ and such that if
$c\geq c_0$ then \whp\ \2G\ finds a 2-matching $M$ with $\k(M)=O(\log n)$.
(This paper gives an analytic proof that $c_0\leq 15$. We have a numerical proof that $c_0\leq 2.5$).
\end{Theorem}
Given this theorem, we will show how we can use this and the result of \cite{FP} to show
\begin{Theorem}\label{th2}\
If $c\geq c_0$ then \whp\ a Hamilton cycle can be found in $O(n^{1.5+o(1)})$ time.
\end{Theorem}

{\bf Acknowledgement:} I would like to thank my colleague Boris Pittel for his
help with this paper. He ought to be a co-author, but he has declined to do so.
\section{Outline of the paper}
As already indicated, the idea is to use a greedy algorithm to find a good 2-matching
and then transform it into a Hamilton cycle.
We will first give an over-view of our greedy algorithm. As we proceed, we select edges
to add to our 2-matching $M$. Thus $M$ consists of paths and cycles
(and isolated vertices). Vertices of the cycles and vertices
interior to the paths get deleted from the current graph, which we
denote by $\G$. No more edges can be added incident
to these interior vertices. Thus the paths can usefully be thought
of as being contracted to the set of edges of a matching
$\M$ on the remaining vertices of $\G$. This matching is not part of $\G$.
We keep track of the vertices covered by $\M$ by using a
0/1 vector $b$ so that for vertex $v$, $b(v)$ is the indicator that $v$ is covered by $\M$.
Thus when $v$ is still included in $\G$ and $b(v)=1$, it will be
the end-point of a path in the current 2-matching $M$.

The greedy algorithm first tries to cover vertices of degree at most two that are not covered
by $M$ or vertices of degree one that are covered by $M$. These choices are forced. When
there are no such vertices, we choose an edge at random. We make sure that one of the end-points
$u,v$ of the chosen edge has $b$-value zero. The aim here is to try to quickly ensure
that $b(v)=1$ for all vertices of $\G$. This will essentially reduce the
problem to that of finding another (near)
perfect matching in $\G$. The first phase of the algorithm finishes when
all of the vertices that remain have $b$-value one. This necessarily means that the contracted
paths form a matching of the
graph $\G$ that remains at this stage. Furthermore, we will see that $\G$ is distributed as
$G_{\n,\m}^{\d\geq 2}$ for some $\n,\m$ and then we construct another (near) perfect
matching $M^{**}$ of $\G$ by using the linear expected time algorithm of \cite{CFM}. We put
$M$ and $M^{**}$ together to create a 2-matching along with the cycles that have been deleted.
Note that some vertices may have become isolated during the construction of $M$ and these
will form single components of our 2-matching. The union of two random (near) perfect matchings
is likely to have $O(\log n)$ components. Full details of this algorithm are given in
Section~\ref{alg}.

Once we have described the algorithm, we can begin its analysis.
We first describe the random graph model that
we will use. We call it the {\em Random Sequence model}. It was first used in Bollob\'as and
Frieze \cite{BoFrM} and independently in
Chvat\'al \cite{Ch}. We used it in \cite{AFP} for our analysis of the
Karp-Sipser algorithm. We prove the truncated Poisson nature of the degree
sequence of the graph $\G$ that remains at each stage in Section
\ref{mod}.
We then, in Section \ref{alg}, give a detailed description of \2G.
In Section \ref{uni} we show that the distribution of the evolving graph $\G$ can be
succinctly
described by a 6-component vector $\bv=(y_1,y_2,z_1,y,z,\m)$ that evolves as a Markov chain.
Here $y_j,j=1,2$ denotes the number of vertices of degree $j$ that are not incident with $M$ and
$z_1$ denotes the number of vertices of degree one that are incident with $M$.
$y$ denotes the number of vertices of degree at least three that are not incident with $M$ and
$z$ denotes the number of vertices of degree at least two that are incident with $M$.
$\m$ denotes the number of edges.
It is important
to keep $\z=y_1+y_2+z_1$ small and \2G\ will attempt to handle such vertices when $\z>0$.
In this way we keep $\z$ small \whp\ throughout the algorithm and this will mean that
the final 2-matching produced will have few components.
Section \ref{cec} first describes the (approximate) transition
probabilities of this chain. There are four types of step in \2G\ that depend on which if any
of $y_1,y_2,z_1$ are positive. Thus there are four sets of transition probabilities.
Given the expected changes in \bv, we first show that in all cases the expected change in $\z$
is negative, when $\z$ is positive. This
indicates that $\z$ will not get large and a high probability
polylog bound is proven.

We are using the differential equation method and Section
\ref{diff} describes the sets of differential equations that can be used to
track the progress of the algorithm \whp. The parameters for these equations will be
$\hat{{\bf v}}=(\hy_1,\hy_2,\hz,\hy,\hz,\hm)$. There are four sets of equations corresponding to the
four types of step in \2G. It is important to know the
proportion of each type of step over a small interval. We thus consider
a {\em sliding trajectory} i.e. a weighted sum of these
four sets of equations. The weights are chosen so that in the
weighted set of equations we have $\hy_1'=\hy_2'=\hz_1'=0$.
This is in line with the fact that $\hy,\hz,\hm\gg\z$ for most of the
algorithm. We verify that the expressions for the weights are
non-negative. We then verify that \whp\ the sliding trajectory and the
process parameters remain close.

Our next aim is to show that \whp\ there is a time $T$ such that
$y(T)=0,z(T)=\Omega(n)$. It would therefore be most natural to show
that for the sliding trajectory, there is a time $\hT$ such that
$\hy(\hT)=0,\hz(\hT)=\Omega(n)$. The equations for the
sliding trajectory are complicated and we have not been able to
do this directly. Instead, we have set up an approximate
system of equations (in parameters $\tiy,\tiz,\tim$) that are
close when $c\geq 15$. We can prove these parameters stay
close to $\hy,\hz,\hm$ and that there is a time $\tiT$ such
that $\tiy(\tiT)=0,\tiz(\tiT)=\Omega(n)$. The existence of
$\hT$ is deduced from this and then we can deduce the existence
of $T$. We then in Section \ref{noc} show that \whp\
\2G\ creates a matching with $O(\log n)$ components, completing the proof of Theorem \ref{th1}.

Section \ref{posa} shows how to use an extension-rotation
procedure on our graph $G$ to find a Hamilton cycle within the claimed
time bounds. This procedure works by extending paths one edge at a time and using an operation
called a rotation to increase the number of chances of extending a path.
It is not guaranteed to extend a path, even if it is possible some other way.
There is the notion of a {\em booster}. This is a non-edge whose addition will allow progress in the
extension-rotation algorithm. The companion paper \cite{FP} shows that
for $c\geq 2.67$ there will \whp\ always be many boosters.
To get the non-edges we first randomly choose
$s=n^{1/2}\log^{-2}n$ random edges $X$ of $G$, none of which
are incident with a vertex of degree three. We then
write $G=G'+X$ and argue in Section \ref{remove} that the
pair $(G',X)$ can be replaced by $(H,Y)$ where
$H=G_{n,cn-s}^{\d\geq 3}$ and $Y$ is a random set of edges
disjoint from $E(H)$. We then argue in Section \ref{exro}
that \whp\ $Y$ contains enough boosters to create a Hamilton cycle within the claimed time bound.

Section \ref{final} contains some concluding remarks.
\section{Random Sequence Model}
\label{mod}
A small change of model will simplify the analysis.
Given a sequence $\bx = x_1,x_2,\ldots,x_{2M}\in [n]^{2M}$ of $2M$ integers between 1 and
$N$ we can define a (multi)-graph
$G_{\bx}=G_\bx(N,M)$ with vertex set $[N]$ and edge set $\{(x_{2i-1},x_{2i}):1\leq
i\leq M\}$. The degree $d_\bx(v)$ of $v\in [N]$ is given by
$$d_\bx(v)=|\set{j\in [2M]:x_j=v}|.$$
If $\bx$ is chosen randomly
from $[N]^{2M}$ then $G_{\bx}$ is close in distribution to $G_{N,M}$. Indeed,
conditional on being simple, $G_{\bx}$ is distributed as $G_{N,M}$. To see this, note that
if $G_{\bx}$ is simple then it has vertex set $[N]$ and $M$ edges.
Also, there are $M!2^M$ distinct equally likely
values of $\bx$ which yield the same graph.

Our situation is complicated by there being lower bounds of $2,3$ respectively
on the minimum degree
in two disjoint sets $J_2,J_3\subseteq [N]$. The vertices in $J_0=[N]\setminus J_2\cup J_3$
are of fixed degree bounded degree and the sum of their degrees is $D=o(N)$.
So we let
$$[N]^{2M}_{J_2,J_3;D}=\\
\{\bx\in [N]^{2M}:d_\bx(j)\geq i\text{ for }j\in J_i,\,i=2,3\text{ and }
\sum_{j\in J_0}d_\bx(j)=D\}.$$
Let $G=G(N,M,J_2,J_3;D)$
be the multi-graph $G_\bx$ for $\bx$ chosen uniformly from
$[N]^{2M}_{J_2,J_3;D}$. It is clear then that conditional on being simple,
$G(n,m,\emptyset,[n];0)$
has the same distribution as $G_{n,m}^{\d\geq 3}$. It is important therefore to estimate the
probability that this graph is simple. For this and other reasons,
we need to have an understanding of
the degree sequence $d_\bx$ when $\bx$ is drawn uniformly from $[N]^{2M}_{J_2,J_3;D}$. Let
$$f_k(\l)=e^\l-\sum_{i=0}^{k-1}\frac{\l^i}{i!}$$
for $k\geq 0$.
\begin{lemma}
\label{lem3}
Let $\bx$ be chosen randomly from
$[N]^{2M}_{J_2,J_3;D}$. For $i=2,3$
let $Z_j\,(j\in [J_i])$
be independent copies of a {\em truncated Poisson} random variable $\cP_i$, where
$$\Pr(\cP_i=t)=\frac{{\l}^t}{t!f_i({\l})},\hspace{1in}t=i,i+1,\ldots\ .$$
Here ${\l}$ satisfies
\begin{equation}\label{2}
\sum_{i=2}^3\frac{{\l}f_{i-1}({\l})}{f_i({\l})}|J_i|=2M-D.
\end{equation}
For $j\in J_0$, $Z_j=d_j$ is a constant and $\sum_{j\in J_0}d_j=D$.
Then $\{d_\bx(j)\}_{j\in [N]}$ is distributed as $\{Z_j\}_{j\in [N]}$ conditional on
$Z=\sum_{j\in [n]}Z_j=2M$.
\end{lemma}
\proofstart
Note first that the value of ${\l}$ in (\ref{2}) is chosen so that
$$\E(Z)=2M.$$
Fix $J_0,J_2,J_3$ and $\boldsymbol{\xi}=(\xi_1,\xi_2,\ldots,\xi_N)$ such that
$\xi_j=d_j$ for $j\in J_0$ and $\xi_j\geq k$ for $k=2,3$.   Then,
$$\Pr(d_{\bx}= \boldsymbol{\xi}) =
\left( \frac{(2M)!}{\xi_1! \xi_2! \ldots \xi_N! }\right)
\bigg/
\left( \sum_{\bx\in [N]^{2M}_{J_2,J_3;D}} \frac{(2M)!}{x_1! x_2! \ldots x_N!} \right).$$
On the other hand,
\begin{align*}
&\Pr \left((Z_1,Z_2,\ldots,Z_N)= \boldsymbol{\xi} \bigg|\; \sum_{j\in[N]} Z_j =  2M\right)=\\
&\left( (2M)!\prod_{j\in J_0}\frac{1}{d_j!}\prod_{i=2}^3\prod_{j\in J_i}
\frac{{\l}^{\xi_j}}{f_i({\l}) \xi_j!}
\right)
\bigg/
\left( \sum_{\bx\in [N]^{2M}_{J_2,J_3;D}} (2M)!\prod_{j\in J_0}\frac{1}{d_j!}
\prod_{i=2}^3\prod_{j\in J_i}
\frac{{\l}^{x_j}}{f_i({\l}) x_j!}
\right)\\
&=\left(   \frac{ \prod_{i=2}^3f_i({\l})^{-|J_i|} {\l}^{2M}}{\xi_1! \xi_2! \ldots \xi_N!} \right)
\bigg/
\left( \sum_{\bx\in [N]^{2M}_{J_2,J_3;D}} \frac{\prod_{i=2}^3f_i({\l})^{-|J_i|} {\l}^{2M}}{x_1! x_2!
\ldots x_N!} \right)\\
&=\Pr(d_\bx= \boldsymbol{\xi}).
\end{align*}
\proofend

To use Lemma \ref{lem3} for the approximation of vertex degrees distributions
we need to have  sharp estimates of the probability that $Z$ is
close to its mean $2M$. In particular we need sharp estimates of
$\Pr(Z=2M)$ and
$\Pr(Z-Z_1=2M-k)$, for $k=o(N)$. These estimates are possible
precisely because $\E(Z)=2M$.
Using the special properties of $Z$, we
can refine a standard argument to show (Appendix 1) that where $N_\ell=|J_\ell|$ and
$\N=N_2+N_3$ and
the variances are
\begin{multline}\label{30}
\s_\ell^2=\frac{f_\ell({\l})({\l}^2f_{\ell-2}({\l})+{\l}f_{\ell-1}({\l}))-{\l}^2f_{\ell-1}({\l})^2}
{f_\ell({\l})^2}={\l}\frac{d}{d{\l}}\bfrac{{\l}f_{\ell-1}({\l})}{f_\ell({\l})}
\\ \text{ and }\s^2=\frac{1}{\N}\sum_{\ell=2}^3N_\ell\s_\ell^2
\end{multline}
that if $\N\s^2\rightarrow \infty$ and $k=O(\sqrt{\N}\s)$
then
\beq{ll1}
\Pr\left(Z=2M-k\right)=\frac{1}{\s\sqrt{2\p \N}}\left(1+
O\bfrac{k^2+1}{\N\s^{2}}\right).
\eeq

A proof for $J_2=[N]$ was given in the appendix of \cite{AFP}.
We need to modify the proof in a trivial way. Given \eqref{ll1}
and
$$\s_\ell^2=O({\l}),\qquad\ell=2,3,$$
we obtain
\begin{lemma}
\label{lem4}
Let $\bx$ be chosen randomly from
$[N]^{2M}_{J_2,J_3;D}$.
\begin{description}
\item[(a)] Assume that $\log \N=O((\N {\l})^{1/2})$. For every
$j\in J_\ell$ and $\ell\leq k\leq \log \N$,
\beq{f1}
\Pr(d_\bx(j)=k)=\frac{{\l}^k}{k!f_\ell({\l})}
\left(1+O\left(\frac{k^2+1}{\N {\l}}\right)\right).
\eeq
Furthermore, for all $\ell_1,\ell_2\in\set{2,3}$
and $j_1\in J_{\ell_1},j_2\in J_{\ell_2},\,j_1\neq j_2$, and
$\ell_i\leq k_i\leq \log \N$,
\beq{f2}
\Pr(d_\bx(j_1)=k_1,d_\bx(j_2)=k_2)=\frac{{\l}^{k_1}}{k_1!f_{\ell_1}({\l})}\frac{{\l}^{k_2}}{k_2!
f_{\ell_2}({\l})}\left(1+O\bfrac{\log^2 \N}{\N {\l}}\right).
\eeq
\item[(b)]
\beq{maxdegree}
d_\bx(j)\leq\frac{\log N}{(\log\log N)^{1/2}} \quad\qs\footnote{An event
$\cE=\cE(\N)$
occurs quite surely
(\qs, in short) if $\Pr(\cE)=1-O(N^{-a})$ for any constant $a>0$}
\eeq
for all $j\in J_2\cup J_3$.
\end{description}
\end{lemma}
\proofstart
Assume that $j=1\notin J_0$. Then
\begin{eqnarray*}
\Pr(d_\bx(1)=k)&=&\frac{\Pr\left(Z_1=k\mbox{ and
}\sum_{i=1}^NZ_i=2M\right)}{\Pr\left(\sum_{i=1}^NZ_i=2M\right)}\\
&=&\frac{{\l}^k}{k!f_\ell({\l})}\frac{\Pr\left(\sum_{i=2}^NZ_i=2M-k\right)}{
\Pr\left(\sum_{i=1}^NZ_i=2M\right)}.
\end{eqnarray*}
Likewise, with $j_1=1,j_2=2$,
$$
\Pr(d_\bx(1)=k_1,d_\bx(2)=k_2)=\frac{{\l}^{k_1}}{k_1!f_{\ell_1}({\l})}\frac{{\l}^{k_2}}{k_2!f_{\ell_2}({\l})}
\frac{\Pr\left(\sum_{i=3}^NZ_i=2M-k_1-k_2\right)}{\Pr\left(\sum_{i=1}^NZ_i=2M\right
)}.\nonumber
$$
Statement (a) follows immediately from \eqref{ll1} and (b) follows from
simple estimations.
\proofend

Let $\n_\bx^\ell(s)$ denote the number of vertices in $J_\ell,\ell=2,3$ of degree
$s$ in $G_\bx$.
Equation (\ref{ll1}) and a standard tail estimate for the binomial
distribution shows
\begin{lemma}
\label{lem4x}
Suppose that $\log \N=O((\N {\l})^{1/2})$ and $N_\ell\to\infty$ with $N$.
Let $\bx$ be chosen randomly from
$[N]^{2M}_{J_2,J_3;D}$. Then \qs,
\beq{degconc}
\cD(\bx)=
\left\{\left|\n_\bx^\ell(j)-\frac{N_\ell {\l}^j}{j!f({\l})}\right|
\leq \brac{1+\bfrac{N_\ell {\l}^j}{j!f({\l})}^{1/2}}\log^2 N,\ k\leq j\leq \log N\right\}.
\eeq
\end{lemma}
\proofend

We can now show $G_\bx$, $\bx\in [n]^{2m}_{\emptyset,[n];0}$ is a good model for
$G_{n,m}^{\d\geq 3}$. For this we only need to show now that
\beq{simpx}
\Pr(G_\bx\text{ is simple})=\Omega(1).
 \eeq
For this we can use a result of McKay \cite{McK}. If we fix the degree sequence
of $\bx$ then $\bx$ itself is just a random permutation of the multi-graph in which each $j\in [n]$
appears $d_\bx(j)$ times. This in fact is another way of looking at the Configuration model of
Bollob\'as \cite{B2}. The reference \cite{McK} shows that the probability $G_\bx$ is simple is
asymptotically equal to $e^{-(1+o(1))\r(\r+1)}$ where $\r=m_2/m$ and
$m_2=\sum_{j\in [n]}d_{\bx}(j)(d_{\bx}(j)-1)$. One consequence of the exponential tails in
Lemma \ref{lem4x} is that $m_2=O(m)$. This implies that $\r=O(1)$ and hence that \eqref{simpx}
holds. We can thus use the Random Sequence Model to prove the occurrence of high probability
events in $G_{n,m}^{\d\geq 3}$.

With this in hand, we can now proceed to describe our 2-matching algorithm.
\section{Greedy Algorithm}\label{alg}
Our algorithm will be applied to the random graph $G=G_{n,m}^{\d\geq 3}$ and analysed in the context
of $G_\bx$. As the algorithm progresses, it makes changes to $G$ and we let $\G$ denote the current
state of $G$. The algorithm grows a 2-matching $M$ and for
$v\in [n]$ we let $b(v)$ be the 0/1 indicator for
vertex $v$ being incident to an edge of $M$. We let
\begin{itemize}
\item $\m$ be the number of edges in $\G$,
\item $V_{0,j}=\set{v\in [n]:d_\G(v)=0,\,b(v)=j}$, $j=0,1$,
\item $Y_k=\set{v\in [n]:d_\G(v)=k\text{ and }b(v)=0}$, $k=1,2$,
\item $Z_1=\set{v\in [n]:d_\G(v)=1\text{ and }b(v)=1}$,
\item $Y=\set{v\in [n]:d_\G(v)\geq 3\text{ and }b(v)=0}$,\qquad
This is $J_3$ of Section \ref{mod}.
\item $Z=\set{v\in [n]:d_\G(v)\geq2\text{ and }b(v)=1}$,\qquad
This is $J_2$ of Section \ref{mod}.
\item $M$ is the set of edges in the current 2-matching.
\item $\M$ is the matching induced by the path components of $M$ i.e. if $P\subseteq M$ is a
path from $x$ to $y$ then $(x,y)$ will be an edge of $\M$ and the internal edges of $P$ will
have been deleted from $\G$.
\end{itemize}
Observe that the sequence $\bb=(b(v))$ is determined by $V_{0,0},V_{0,1},Y_1,Y_2,Z_1,Y,Z$.

If $Y_1\neq 0$ then we choose $v\in Y_1$ and add
the edge incident to $v$ to $M$, because
doing so is not a mistake i.e. there is a maximum size 2-matching of $\G$ that
contains this edge.
If $Y_1=\emptyset$ and $Y_2\neq \emptyset$ then we choose $v\in Y_2$ and add
the two edges incident to $v$ to $M$, because
doing so is also not a mistake i.e. there is a maximum size 2-matching of $\G$ that
contains these edges.
Similarly, if $Y_2=\emptyset$ and $Z_1\neq \emptyset$ we choose $v\in Z_1$ and add the unique
edge of $\G$
incident to $v$ to $M$. When we add an edge to $M$ it can cause vertices of $\G$ to become internal
vertices of paths of $M$ and be deleted from $\G$. In particular, this happens to $v\in Z_1$
in the case just described. When $Y_1=Y_2=Z_1=\emptyset\neq Y$ we choose a random edge incident to
a vertex of $Y$. In this way we hope to end up in a situation where $Y_2=Z_1=Y=\emptyset$
and $|Z|=\Omega(n)$. This has advantages that will be explained later in Section \ref{noc} and we have
only managed to prove that this happens \whp\ when $c\geq 15$. When $Y_1=Y_2=Z_1=Y=\emptyset$ we are looking for
a maximum matching in the graph $\G$ that remains and we can use the results of \cite{CFM}.

We now give details of the steps of \vspace{.1in}

\noindent
{\bf Algorithm \2G:}\vspace{-.1in}
\begin{description}
\item[Step 1(a) $Y_1\neq\emptyset$]\ \\
Choose a random vertex $v$ from $Y_1$. Suppose that its neighbour in $\G$ is
$w$. We add $(v,w)$ to $M$ and move $v$ to $V_{0,1}$.
\begin{enumerate}[(i)]
\item If $b(w)=0$ then we add $(v,w)$ to  $\M$. If $w$ is currently
in $Y$ then move it to $Z$.
If it is currently in $Y_1$ then move it to $V_{0,1}$.
If it is currently in $Y_2$ then move it to $Z_1$. Call this re-assigning $w$.
\item If $b(w)=1$ let $u$ be the other end point of the path $P$ of $M$
that contains $w$. We remove
$(w,u)$ from $\M$ and replace it with $(v,u)$. We move $w$ to $V_{0,1}$ and
make the requisite changes
due to the loss of other edges incident with $w$. Call this {\em tidying up}.
\end{enumerate}

\item[Step 1(b): $Y_1=\emptyset$ and $Y_2\neq\emptyset$]\ \\
Choose a random vertex $v$ from $Y_2$. Suppose that its neighbours in $\G$ are
$w_1,w_2$.

If $w_1=w_2=v$ then we simply delete $v$ from $\G$.
(We are dealing with loops because we are
analysing the algorithm within the context of $G_\bx$.
This case is of course unnecessary when the input is simple i.e.
for $G_{n,m}^{\d\geq k}$).

Continuing with the most likely case, we move $v$ to $V_{0,1}$.
We delete the edges\\ $(v,w_1),(v,w_2)$
from $\G$ and place them into $M$.
In addition,
\begin{enumerate}[(i)]
\item If $b(w_1)=b(w_2)=0$ then we add $(w_1,w_2)$ to $\M$ and put
$b(w_1)=b(w_2)=1$. Re-assign $w_1,w_2$.
\item If $b(w_1)=b(w_2)=1$ let $u_i,i=1,2$ be the other end points of the
paths $P_1,P_2$ of $M$ that contain $w_1,w_2$
respectively.
There are now two possibilities:
\begin{enumerate}[(1)]
\item $u_1=w_2$ and $u_2=w_1$. In this case, adding the two edges
creates a cycle $C=(v,w_1,P_1,w_2,v)$
and we delete the edge $(w_1,w_2)$ from $\M$. Vertices $w_1,w_2$ are deleted from $\G$.
The rest of $C$ has already been deleted. Tidy up.
\item $u_1\neq w_2$ and $u_2\neq w_1$. Adding the two edges creates a path\\
$(u_1,P_1',w_1,v,w_2,P_2,u_2)$ to $M$, where $P_1'$ is the reversal
of $P_1$. We delete the edges $(w_1,u_1),(w_2,u_2)$ from $\M$
and add $(u_1,u_2)$ in their place. Vertices $w_1,w_2$ are deleted from $\G$.  Tidy up.
\end{enumerate}
\item If $b(w_1)=0$ and $b(w_2)=1$ let $u_2$ be the other end point
of the path $P_2$ of $M$ that contains $w_2$. We delete
$(w_2,u_2)$ from $\M$ and replace it with $(w_1,u_2)$. We put $b(w_1)=1$
and re-assign it and delete vertex $w_2$ from $\G$.  Tidy up.
\end{enumerate}

\item[Step 1(c): $Y_2=\emptyset$ and $Z_1\neq\emptyset$]\ \\
Choose a random vertex $v$ from $Z_1$. Let $u$ be the other endpoint of the path $P$ of $M$
that contains $v$. Let $w$ be the unique neighbour of $v$ in $\G$. We delete $v$ from $\G$ and
add the edge $(v,w)$ to $M$. In addition there
are two cases.
\begin{enumerate}[(1)]
\item If $b(w)=0$ then we delete $(v,u)$ from $\M$ and replace it
with $(w,u)$ and put $b(w)=1$ and re-assign $w$.
\item If $b(w)=1$ then let $u$ be the other end-point of the path
containing $w$ in $M$. If $u\neq v$ then we delete vertex $w$ and the edge $(u,w)$ from $\M$ and replace it with
$(u,v)$.  Tidy up. If $u=v$ then we have created a cycle $C$ and we delete it from $\G$ as in Step 1(b)(i)(1).
\end{enumerate}

\item[Step 2: $Y_1=Y_2=Z_1=\emptyset$ and $Y\neq\emptyset$]\ \\
Choose a random edge $(v,w)$ incident with a vertex $v\in Y$.
We delete the edge $(v,w)$ from $\G$ and add it to $M$. We put $b(v)=1$
and move it from $Y$ to $Z$. There are two cases.
\begin{enumerate}[(i)]
\item If $b(w)=0$ then put $b(w)=1$ and move it from $Y$ to $Z$. We
add the edge $(v,w)$ to $\M$.
\item If $b(w)=1$ let $u$ be the other end point of the path in $M$
containing $w$. We delete vertex $w$ and the edge $(u,w)$ from $\M$ and replace it with
$(u,v)$.  Tidy up.
\end{enumerate}
\item[Step 3: $Y_1=Y_2=Z_1=Y=\emptyset$]\ \\
At this point $\G$ will be seen to be distributed as $G_{\n,\m}^{\d\geq 2}$
for some $\n,\m$ where $\m=O(\n)$.
As such, it contains a (near) perfect matching $M^{**}$ \cite{FP} and it can
be found in $O(\n)$ expected time \cite{CFM}.
\end{description}
The output of \2G\ is set of edges in $M\cup M^{**}$.

No explicit mention has been made of vertices contributing to $V_{0,0}$. When we we tidy
up after removing a vertex $w$, any vertex whose sole neighbour is $w$ will be placed in
$V_{0,0}$.

\section{Uniformity}\label{uni}
In the previous section, we described the action of the algorithm as applied to $\G$.
In order to prove a uniformity property, it is as well to consider the changes induced
by the algorithm in terms of \bx.

When an edge is removed
we will replace it in $\bx$ by a pair of $\star$'s. This goes for all of
the edges removed at
an iteration, not just the edges of the 2-matching $M$. Thus at the end of
this and subsequent iterations we will have a sequence in
$\La=([n]\cup \{\star\})^{2m}$ where for all $i$, $x_{2i-1}=\star$ if and
only if $x_{2i}=\star$.
We call such sequences {\em proper}.

We use the same notation as in Section \ref{mod}. Let
$S=S(\bx)=\{i:z_{2i-1}=z_{2i}=\star\}$. Note that the number of edges $\m$ in $G_\bx$
is given by
$$\m=m-|S|.$$
For a
tuple $\bv=(V_{0,0},V_{0,1},Y_1,Y_2,Z_1,Y,Z,S)$
we let $\La_{\bv}$ denote the set of pairs $(\bx,\bb)$ where $\bx\in \La$ is proper and
\begin{itemize}
\item $V_{0,j}=\set{v\in [n]:d_\bx(v)=0,\,b(v)=j}$, $j=0,1$,
\item $Y_k=\set{v\in [n]:d_\bx(v)=k\text{ and }b(v)=0}$, $k=1,2$,
\item $Z_1=\set{v\in [n]:d_\bx(v)=1\text{ and }b(v)=1}$,
\item $Y=\set{v\in [n]:d_\bx(v)\geq 3\text{ and }b(v)=0}$,
\item $Z=\set{v\in [n]:d_\bx(v)\geq2\text{ and }b(v)=1}$.
\item $S=S(\bx)$.
\end{itemize}
(Re-call that \bb\ is determined by \bv).

For vectors $\bx,\bb$ we define $\bv(\bx,\bb)$ by $(\bx,\bb)\in \La_{\bv(\bx,\bb)}$. We also use the
notation $\bx\in \La_{\bv(\bx)}$
when the second component \bb\ is assumed.

Given two sequences $\bx,\bx'\in \La$, we say that $\bx'\subseteq \bx$ if $x_j=\star$
implies $x_j'=\star$.
In which case we define $\by=\bx-\bx'$ by
$$y_j=\begin{cases}
       x_j&\text{ If }x_j\neq \star=x_j'\\ \star&\text{ Otherwise}
      \end{cases}
$$
Thus $\by$ records the changes in going from \bx\ to $\bx'$.

Given two sequences $\bx,\bx'\in \La$ we say that $\bx,\bx'$ are {\em disjoint}
if $x_j\neq\star$ implies that $x_j'=\star$. In which case we define $\by=\bx+\bx'$ by
$$y_j=\begin{cases}
      x_j&\text{ If }x_j\neq \star\\
      x_j'&\text{ If }x_j'\neq \star\\
      \star&\text{ Otherwise}
      \end{cases}
$$
Thus,
\beq{1-1}
\mbox{if $\bx'\subseteq \bx$ then $\bx'$ and $\bx-\bx'$ are disjoint and $\bx=\bx'+(\bx-\bx')$.}
\eeq
Suppose now that $(\bx(0),\bb(0)),(\bx(1),\bb(1)),\ldots,(\bx(t),\bb(t))$ is
the sequence of pairs representing the graphs constructed by the algorithm \2G.
Here $\bx(i-1)\supseteq \bx(i)$ for $i\geq 1$ and so we can define $\by(i)=\bx(i-1)-\bx(i)$.
Suppose that $\bv(i)=\bv(\bx(i))$ for $1\leq i\leq t$
where $\bv(0)=(\emptyset,\emptyset,\emptyset,\emptyset,[n],\emptyset,\emptyset)$
and $\bb(0)=0$.

Let
$$\La_{\bv\mid \bb}=\set{\bx:(\bx,\bb)\in \La_\bv}.$$
\begin{lemma}
\label{lem2}
Suppose that $\bx(0)$ is a random member of $\La_{\bv(0)\mid\bb(0)}$. Then
given\\
 $\bv(0),\bv(1),\ldots,\bv(t)$,
the vector $\bx(t)$ is a random member of
$\La_{\bv(t)\mid\bb(t)}$ for all
$t\geq 0$, that is, the distribution of $\bx(t)$ is uniform,
conditional on the edges deleted in the first $t$ steps. (Note
that $\bb(t)$ is fixed by $\bv(t)$ here).
\end{lemma}
\proofstart
We prove this by induction on $t$. It is trivially true for $t=0$. Fix
$t\geq 0,\bx(t),\bb(t),\bx(t+1),\bb(t+1)$.
We define a sequence $\bx(t)=\bz_1,\bz_2,\ldots,\bz_s=\bx(t+1)$
where $\bz_{i+1}$ is obtained from $\bz_i$ by a {\em basic step}

\vspace{.2in}
\noindent
{\bf Basic Step}: Given $\bx,\bb$ and $\bv=\bv(\bx,\bb)$ we create
new sequences $\bx'=A_j(\bx),\bb'=B_j(\bb)$ and $\bv'=\bv(\bx',\bb')$.
Let $\bw=\bx-\bx'$.
A basic step corresponds to replacing the edge $(w_{2j-1},w_{2j})$
by an edge of the matching $M$, for some index $j$.
Let $u=w_{2j-1},v=w_{2j}$.
\begin{description}
\item[Case 1:] Here
we assume $b(u)=b(v)=0$.\\
Replace $x_{2j-1},x_{2j}$ by $\star$'s and put $b(u)=b(v)=1$.
\item[Case 2:] Here
we assume $b(u)=0,\,b(v)=1$.\\
Replace $x_{2k-1},x_{2k}$ by $\star$'s for every $k$ such that $v\in \set{x_{2k-1},x_{2k}}$
and put $b(u)=1$.
\item[Case 3:] Here
we assume $b(u)=b(v)=1$.\\
Replace $w_{2k-1},w_{2k}$ by $\star$'s for every $k$ such that
$\set{u,v}\cap \set{w_{2k-1},w_{2k}}
\neq\emptyset$.
\end{description}
\begin{claim}\label{cl11}
Suppose that $\bx'=A_j(\bx)$ and $\by=\bx-\bx'$ and $\bb'=B_j(\bb)$.
Then the map $\f:\bz\in \La_{\bv(\bx,\bb)}^{\by}\to(\bz-\by,\bb')$ is 1-1 and
each $(\bz',\bb')\in \La_{\bv(\bx',\bb')}$ is the image under $\f$
of a unique member of $\La_{\bv(\bx,\bb)}^\by$,
where $\La_{\bv(\bx,\bb)}^\by=\set{(\bz,\bb)\in \La_{\bv(\bx,\bb)}:\;\bz\supseteq\by}$.
\end{claim}
{\bf Proof of Claim \ref{cl11}}. Equation \eqref{1-1} implies that $\f$ is 1-1.
Let $\bv=\bv(\bx,\bb)$ and $\bv'=\bv(\bx',\bb')$.
Choose $(\bw,\bb')\in \La_{\bv'}$. Because $S'$ is determined by
$\bv'$, we see that \by\ and \bw\ are necessarily disjoint
and we simply have to check that if $\sx=\bw+\by$ then $(\sx,\bb)\in \La_{\bv}$.
But in all cases, $\bv(\sx,\bb)$ is determined by $\bv'$ and $\by$ and
this implies that $\bv(\sx,\bb)=\bv(\bx,\bb)$.

This statement is the crux of the proof and we should perhaps justify it a little more.
Suppose then that we are given $\bv'$ (and hence $\bb'$) and
$\by$ and \bb. Observe that this determines
$d_{\bx^*}(v)$ for all $v\in V_{0,0}'\cup V_{0,1}'\cup Y_1'\cup Y_2'\cup Z_1'$. Together
with $b(v)$ this determines the place of $v$ in the partition defined by \bv. Now
$Y'\subseteq Y$ and it only remains to deal with $v\in Z'$. If $d_\by(v)>0$ then $v\in
Y\cup Z$ and $b(v)$ determines which of the sets $v$ is in. If $d_\by(v)=0$ and $b(v)=1$
then $v\in Z$. If $d_\by(v)=0$ and $b(v)=0$ then $v\in Y$. This is because $b(v)=0$ and
$b'(v)=1$ implies that we have put one of the edges incident with $v$ into $M$.
\\
{\bf End of proof of Claim \ref{cl11}}

The claim implies (inductively) that if $\bx$ is a uniform random member
of $\La_{\bv\mid\bb}$ and we do a sequence of basic
steps involving the ``deletion'' of $\by_1,\by_2,\ldots,\by_s$ where
$\by_{i+1}\subseteq \bx-\by_1-\cdots
\by_i$, then $\bx'=\bx-\by_1-\cdots -\by_s$ is a uniform
random member of $\La_{\bv'\mid\bb'}$, where $\bv'=\bv(\bx',\bb')$
for some $\bb'$. This will imply Lemma \ref{lem2} once
we check that a step of \2G\ can be broken into basic steps.

First consider Step 1(a). First we choose a vertex in $x\in Y_1$.
Then we apply Case 1 or 2 with probabilities determined by \bv.

Now consider Step 1(b).
First we choose a vertex in
$x\in Y_2$. We can then replace one of the
edges incident with $x$ by a matching edge. We apply Case 1 or
Case 2 with probabilities determined by \bv.
After this we apply Case 2 or Case 3 with probabilities determined by \bv.

For Step 1(c)
we apply one of Case 2 or Case 3 with probabilities determined by \bv.

For Step 2, we apply one of Case 1 or Case 2 with probabilities determined by \bv.

This completes the proof of Lemma \ref{lem2}.
\proofend

As a consequence
\begin{lemma}
\label{lem3a}
The random sequence $\bv(t),\,t=0,1,2,\ldots,$ is a Markov chain.
\end{lemma}
\proofstart
Slightly abusing notation,
\begin{align*}
&\Pr(\bv(t+1)\mid \bv(0),\ldots,\bv(t))\\
&=\sum_{\bw'\in
\La_{\bv(t+1)}}\Pr(\bw'\mid
\bv(0),\ldots,\bv(t))\\
&=\sum_{\bw'\in \La_{\bv(t+1)}}\sum_{\bw\in \La_{\bv(t)}}\Pr(\bw',\bw\mid
\bv(0),\ldots,\bv(t))\\
&=\sum_{\bx'\in \La_{\bv(t+1)}}\sum_{\bw\in \La_{\bv(t)}}
\Pr(\bw'\mid \bv(0),\ldots,\bv(t-1),\bw)
\Pr(\bw\mid\bv(0),\ldots,\bv(t))\\
&=\sum_{\bw'\in \La_{\bv(t+1)}}\sum_{\bw\in \La_{\bv(t)}}\Pr(\bw'\mid
\bw)|\La_{\bv(t)}|^{-1},\qquad using\ Lemma\ \ref{lem2}.
\end{align*}
which depends only on $\bv(t),\bv(t+1)$.
\proofend

We now let
$$|\bv|=\set{|V_{0,0}|,|V_{0,1}|,|Y_1|,|Y_2|,|Z_1|,|Y|,|Z|,|S|}.$$
Then we let $\La_{|\bv|}$ denote the set of $(\bx,\bb)\in \La$ with
$|\bv(\bx,\bb)|=|\bv|$ and we let\\ $\La_{|\bv|\,\mid\bb}=\set{\bx:(\bx,\bb)\in \La_{|\bv|}}$.

It then follows from Lemma \ref{lem3a} that by symmetry,
\begin{lemma}
\label{|lem2|}
The random sequence $|\bv(t)|,\,t=0,1,2,\ldots,$ is a Markov chain.
\end{lemma}

A component of a graph
is {\em trivial} if it consists of a single isolated vertex.

\begin{lemma}\label{few}
{\bf Whp} the number of non-trivial components of the graph induced by $M\cup M^{**}$ is $O(\log n)$.
\end{lemma}
\proofstart
Lemma 3 of Frieze and {\L}uczak \cite{FL} proves that \whp\ the union of
two random (near) perfect matchings of $[n]$ has at most
$3\log n$ components. Lemma \ref{lem2} implies that at the end of Phase 1,
$\G$ is a copy of $G_{\n,\m}^{\d\geq 2}$, independent of $\M$. In which case the (near) perfect
matching of $\G$ is independent of $\M$ and we can apply \cite{FL}.
\proofend
\section{Conditional expected changes}\label{cec}
We now set up a system of differential equations that
closely describe the path taken by the parameters of
Algorithm \2G, as applied to $G_\bx$ where $\bx$
is chosen randomly from $[n]^{2m}_{\emptyset,[n];0}$.
We introduce the following notation: At some point in the algorithm, the state of
$\G$ is described by $\bx\in [n]^{2M}_{J_2,J_3;D}$, together with an indicator vector \bb. We let
$y_i=\card{\set{v:d_\bx(v)=i\text{ and }b(v)=0}}$ and
let $z_i=\card{\set{v:d_\bx(v)=i\text{ and }b(v)=1}}$ for $i\ge0$.
We let $y=\sum_{i\ge 3}y_i$ and $z=\sum_{i\ge 2}z_i$ and let
$2\m=\sum_{i\geq 0}i(y_i+z_i)$ be the total degree. Thus in
the notation of Section \ref{alg} we have $y_i=|Y_i|, i=1,2,\,J_3=Y,N_3=y,\,
z_1=|Z_1|,\,J_2=Z,N_2=z,D=y_1+2y_2+z_1, M=\m$. Then it follows
from Lemma \ref{lem4x}, that as long as $(y+z)\l=\Omega(\log^2n)$, we have \qs,
\begin{equation}\label{1}
y_k \approx \frac{\la^k}{k!f_3(\la)},\,\, (k\ge 3);\quad
z_k\approx\,\frac{\la^k}{k!f_{2}(\la)},\,\,(k\ge 2).
\end{equation}
Here $\la$ is the root of
\begin{equation}\label{3}
y\,\frac{\la f_2(\la)}{f_3(\la)}+z\,\frac{\la f_1(\la)}{f_2(\la)}=2\mu-y_1-2y_2-z_1.
\end{equation}
\proofend

{\bf Notational Convention:} There are a large number of parameters that change as\\ \2G\
progresses. Our convention will be that if we write a parameter $\xi$ then by default
it means $\xi(t)$, the value of $\xi$ after $t$ steps of the algorithm. Thus the initial
value of $\xi$ will be $\xi(0)$. When $\xi$ is evaluated at a different point, we make this
explicit.

We now keep track of the expected changes in $\bv=(y_1,y_2,y,z_1,z_2,\m)$ due to one step of \2G.
These expectations are conditional on the current values of $\bb$ and the degree sequence \bd.
 We let $N=y+z$, which is a small departure from the notation of Section \ref{mod}.
In the following sequence of equations, $\xi'=\xi(t+1)$ represents the value of parameter $\xi$
after the corresponding step of \2G.
\begin{lemma}\label{lemchanges}
The following are the expected one step changes in the parameters\\
 $(y_1,y_2,y,z_1,z,\m)$. We will compute them conditional on the degree sequence \bd\
 and on $|\bv|$. We give both, because the first are more transparent and the second
 are what is needed. The error terms $\e_?$ are the consequence of multi-edges and we will
 argue that they are small. We take
$$N=y+z.$$

{\bf  Step 1.\/} $y_1+y_2+z_1>0$.
\si
{\bf Step 1(a)\/}. $y_1>0$.
\begin{eqnarray}
\ex[y_1^\prime-y_1\mid \bb,\bd]&=&-1-\brac{\frac{y_1}{2\mu}+
\sum_{k\ge 2}\frac{kz_k}{2\mu}\,(k-1)\frac{y_1}{2\mu}}
+\sum_{k\ge 2}\frac{kz_k}{2\mu}\,(k-1)\frac{2y_2}{2\mu}+\e_{\ref{04x}}.\label{04x}\\
\ex[y_1^\prime-y_1\mid
|\bv|]&=&-1-\frac{y_1}{2\mu}-\frac{y_1z}{4\mu^2}\frac{\la^2f_0(\la)}{f_2(\la)}
+\frac{y_2z}{2\mu^2}\,\frac{\la^2f_0(\la)}{f_2(\la)}+O\bfrac{\log^2N}{\l N}\label{04xq}\\
\nonumber\\ \nonumber\\
\ex[y_2^\prime-y_2\mid \bb,\bd]&=&-\brac{\frac{2y_2}{2\mu}+
\sum_{k\ge 2}\frac{kz_k}{2\mu}\,(k-1)\frac{2y_2}{2\mu}}
+\sum_{k\ge 2}\frac{kz_k}{2\mu}\,(k-1)\frac{3y_3}{2\mu}+\e_{\ref{04}}.\label{04}\\
\ex[y_2^\prime-y_2\mid |\bv|]&=&-\frac{y_2}{\mu}-
\frac{y_2z}{2\mu^2}\frac{\la^2f_0(\la)}{f_2(\la)}
+\frac{yz}{8\mu^2}\frac{\la^3}{f_3(\la)}\,\frac{\la^2f_0(\la)}{f_2(\la)}
+O\bfrac{\log^2N}{\l N}\label{04q}\\
\nonumber\\ \nonumber\\
\ex[z_1^\prime-z_1\mid \bb,\bd]&=&-\brac{\frac{z_1}{2\mu}+
\sum_{k\ge 2}\frac{kz_k}{2\mu}\,(k-1)\frac{z_1}{2\mu}}
+\sum_{k\ge 2}\frac{kz_k}{2\mu}\,(k-1)\frac{2z_2}{2\mu}+\e_{\ref{05}}.\label{05}\\
\ex[z_1^\prime-z_1\mid |\bv|]&=&-\frac{z_1}{2\mu}
-\frac{z_1z}{4\mu^2}\,\frac{\la^2f_0(\la)}{f_2(\la)}
+\frac{z^2}{4\mu^2}\frac{\la^4f_0(\la)}{f_2(\la)^2}+O\bfrac{\log^2N}{\l N}.
\label{05q}\\
\nonumber\\ \nonumber\\
\ex[y^\prime-y\mid \bb,\bd]&=&-\brac{\sum_{k\ge 3}\frac{ky_k}{2\mu}+
\sum_{k\ge 2}\frac{kz_k}{2\mu}\,(k-1)\frac{3y_3}{2\mu}}+\e_{\ref{06}}.\label{06}\\
\ex[y^\prime-y\mid |\bv|]&=&-\frac{y}{2\mu}\,\frac{\la
 f_2(\la)}{f_3(\la)}-\frac{yz}{8\mu^2}\,\frac{\la^3}{f_3(\la)}\,
\,\frac{\la^2f_0(\la)}{f_2(\la)}+O\bfrac{\log^2N}{\l N}.
\label{06q}\\
\nonumber\\ \nonumber\\
\ex[z^\prime-z\mid \bb,\bd]&=&\sum_{k\geq 3}\frac{ky_k}{2\mu}-\sum_{k\ge 2}
\frac{kz_k}{2\mu}-\sum_{k\ge 2}\frac{kz_k}{2\mu}\,(k-1)\frac{2z_2}{2\mu}
+\e_{\ref{07}}.\label{07}\\
\ex[z^\prime-z\mid |\bv|]&=&\frac{y}{2\mu}
\frac{\la f_2(\la)}{f_3(\la)}-\frac{z}{2\mu}\frac{\la
 f_1(\la)}{f_2(\la)}-\frac{z^2}{4\mu^2}\,\frac{\la^4f_0(\la)}{f_2(\la)^2}
 +O\bfrac{\log^2N}{\l N}.\label{07q}\\
\nonumber\\ \nonumber\\
\ex[\mu^\prime-\mu\mid \bb,\bd]&=&-1-\sum_{k\ge
 2}\frac{kz_k}{2\mu}(k-1)+\e_{\ref{08a}}.\label{08a}\\
\ex[\mu^\prime-\mu\mid|\bv|]&=&-1-\frac{z}{2\mu}\frac{\la^2f_0(\la)}{f_2(\la)}
+O\bfrac{\log^2N}{\l N}.\label{08q}
\end{eqnarray}
{\bf Step 1(b)\/}. $y_1=0,y_2>0$.
\begin{eqnarray}
\ex[y_1^\prime-y_1\mid \bb,\bd]&=&
2\sum_{k\ge 2}\frac{kz_k}{2\mu}\,(k-1)\frac{2y_2}{2\mu}+\e_{\ref{4x}}.\label{4x}\\
\ex[y_1^\prime-y_1\mid
 |\bv|]&=&\frac{y_2z}{\mu^2}\,\frac{\la^2f_0(\la)}{f_2(\la)}
+O\bfrac{\log^2N}{\l N}.\label{4xq}\\
\nonumber\\ \nonumber\\
\ex[y_2^\prime-y_2\mid \bb,\bd]&=&-1-2\brac{\frac{2y_2}{2\mu}+
\sum_{k\ge 2}\frac{kz_k}{2\mu}\,(k-1)\frac{2y_2}{2\mu}}
+2\sum_{k\ge 2}\frac{kz_k}{2\mu}\,(k-1)\frac{3y_3}{2\mu}+\e_{\ref{4}}.\label{4}\\
\ex[y_2^\prime-y_2\mid
 |\bv|]&=&-1-\frac{2y_2}{\mu}-\frac{y_2z}{\mu^2}\frac{\la^2f_0(\la)}{f_2(\la)}
+\frac{yz}{4\mu^2}\frac{\la^3}{f_3(\la)}\,\frac{\la^2f_0(\la)}{f_2(\la)}
+O\bfrac{\log^2N}{\l N}.\label{4q}\\
\nonumber\\ \nonumber\\
\ex[z_1^\prime-z_1\mid \bb,\bd]&=&-2\brac{\frac{z_1}{2\mu}+
\sum_{k\ge 2}\frac{kz_k}{2\mu}\,(k-1)\frac{z_1}{2\mu}}
+2\sum_{k\ge 2}\frac{kz_k}{2\mu}\,(k-1)\frac{2z_2}{2\mu}+\e_{\ref{5}}.\label{5}\\
\ex[z_1^\prime-z_1\mid
|\bv|]&=&-\frac{z_1}{\mu}-\frac{z_1z}{2\mu^2}\,\frac{\la^2f_0(\la)}{f_2(\la)}
+\frac{z^2}{2\mu^2}\frac{\la^4f_0(\la)}{f_2(\la)^2}
+O\bfrac{\log^2N}{\l N}.\label{5q}\\
\nonumber\\ \nonumber\\
\ex[y^\prime-y\mid \bb,\bd]&=&-2\brac{\sum_{k\ge 3}\frac{ky_k}{2\mu}+
\sum_{k\ge 2}\frac{kz_k}{2\mu}\,(k-1)\frac{3y_3}{2\mu}}+\e_{\ref{6}}.\label{6}\\
\ex[y^\prime-y\mid |\bv|]&=&-\frac{y}{\mu}\,\frac{\la
 f_2(\la)}{f_3(\la)}-\frac{yz}{4\mu^2}\,\frac{\la^3}{f_3(\la)}\,
\,\frac{\la^2f_0(\la)}{f_2(\la)}
+O\bfrac{\log^2N}{\l N}.\label{6q}\\
\nonumber\\ \nonumber\\
\ex[z^\prime-z\mid \bb,\bd]&=&2\brac{\sum_{k\geq 3}\frac{ky_k}{2\mu}-
\sum_{k\ge 2}\frac{kz_k}{2\mu}-
\sum_{k\ge 2}\frac{kz_k}{2\mu}\,(k-1)\frac{2z_2}{2\mu}}+\e_{\ref{7}}.\label{7}\\
\ex[z^\prime-z\mid |\bv|]&=&\frac{y}{\mu}
\frac{\la f_2(\la)}{f_3(\la)}-\frac{z}{\mu}\frac{\la
 f_1(\la)}{f_2(\la)}-\frac{z^2}{2\mu^2}\,\frac{\la^4f_0(\la)}{f_2(\la)^2}
+O\bfrac{\log^2N}{\l N}.\label{7q}\\
\nonumber\\ \nonumber\\
\ex[\mu^\prime-\mu\mid \bb,\bd]&=&-2-2\sum_{k\ge
 2}\frac{kz_k}{2\mu}(k-1)+\e_{\ref{8a}}.\label{8a}\\
\ex[\mu^\prime-\mu\mid |\bv|]&=&-2-\frac{z}{\mu}\frac{\la^2f_0(\la)}{f_2(\la)}
+O\bfrac{\log^2N}{\l N}.\label{8aq}\\
\end{eqnarray}
{\bf  Step 1(c).\/} $y_1=y_2=0,z_1>0$.
\begin{eqnarray}
\ex[y_1^\prime-y_1\mid \bb,\bd]&=&O\bfrac{1}{N}.\label{9x}\\
\ex[y_1^\prime-y_1\mid |\bv|]&=&O\bfrac{1}{N}.\label{9xq}\\
\nonumber\\ \nonumber\\
\ex[y_2^\prime-y_2\mid \bb,\bd]&=&\sum_{k\ge 2}\frac{kz_k}{2\mu}\,(k-1)\frac{3y_3}{2\mu}
+\e_{\ref{9}}.\label{9}\\
\ex[y_2^\prime-y_2\mid |\bv|]&=&{\frac{yz}{8\mu^2}\,\frac{\la^3}{f_3(\la)}\,
\frac{\la^2f_0(\la)}{f_2(\la)}}
+O\bfrac{\log^2N}{\l N}.\label{9q}\\
\nonumber\\ \nonumber\\
\ex[z_1^\prime-z_1\mid \bb,\bd]&=&-1-\frac{z_1}{2\mu}-
\sum_{k\ge 2}\frac{kz_k}{2\mu}\,(k-1)\frac{z_1}{2\mu}
+\sum_{k\ge 2}\frac{kz_k}{2\mu}\,(k-1)\frac{2z_2}{2\mu}+\e_{\ref{10}}.\label{10}\\
\ex[z_1^\prime-z_1\mid |\bv|]&=&-1-\frac{z_1}{2\mu}-
\frac{z_1z}{4\mu^2}\,\frac{\la^2f_0(\la)}{f_2(\la)}
+\frac{z^2}{4\mu^2}\frac{\la^4f_0(\la)}{f_2(\la)^2}
+O\bfrac{\log^2N}{\l N}.\label{10q}\\
\nonumber\\ \nonumber\\
\ex[y^\prime-y\mid \bb,\bd]&=&-\sum_{k\ge 3}\frac{ky_k}{2\mu}-
\sum_{k\ge 2}\frac{kz_k}{2\mu}\,(k-1)\frac{3y_3}{2\mu}+\e_{\ref{6a}}.\label{6a}\\
\ex[y^\prime-y\mid |\bv|]&=&-\frac{y}{2\mu}\,
\frac{\la f_2(\la)}{f_3(\la)}-\frac{yz}{8\mu^2}\,\frac{\la^3}{f_3(\la)}\,
\,\frac{\la^2f_0(\la)}{f_2(\la)}
+O\bfrac{\log^2N}{\l N}.\label{6aq}\\
\nonumber\\ \nonumber\\
\ex[z^\prime-z\mid \bb,\bd]&=&\sum_{k\geq 3}\frac{ky_k}{2\mu}-\sum_{k\ge 2}\frac{kz_k}{2\mu}-
\sum_{k\ge 2}\frac{kz_k}{2\mu}\,(k-1)\frac{2z_2}{2\mu}+\e_{\ref{7a}}.\label{7a}\\
\ex[z^\prime-z\mid |\bv|]&=&\frac{y}{2\mu}\,\frac{\la f_2(\la)}{f_3(\la)}
-\frac{z}{2\mu}\frac{\la f_1(\la)}{f_2(\la)}
-\frac{z^2}{4\mu^2}\,\frac{\la^4f_0(\la)}{f_2(\la)^2}
+O\bfrac{\log^2N}{\l N}.\label{7aq}\\
\nonumber\\ \nonumber\\
\ex[\mu^\prime-\mu\mid \bb,\bd]&=&-1-\sum_{k\ge
 2}\frac{kz_k}{2\mu}(k-1)+\e_{\ref{8aa}}.\label{8aa}\\
\ex[\mu^\prime-\mu\mid |\bv|]&=&-1-\frac{z}{2\mu}\frac{\la^2f_0(\la)}{f_2(\la)}
+O\bfrac{\log^2N}{\l N}.\label{8aaq}
\end{eqnarray}

{\bf Step 2.\/} $y_1=y_2=z_1=0$.
\begin{eqnarray}
\ex[y_1^\prime-y_1\mid \bb,\bd]&=&O\bfrac{1}{N}.\label{11x}\\
\ex[y_1^\prime-y_1\mid |\bv|]&=&O\bfrac{1}{N}.\label{11xq}\\
\nonumber\\ \nonumber\\
\ex[y_2^\prime-y_2\mid \bb,\bd]&=&\sum_{k\ge 2}\frac{kz_k}{2\mu}\,(k-1)\frac{3y_3}{2\mu}+
\e_{\ref{11}}.\label{11}\\
\ex[y_2^\prime-y_2\mid |\bv|]&=&\frac{yz}{8\mu^2}\frac{\la^3}{f_3(\la)}
\frac{\la^2f_0(\la)}{f_2(\la)}+O\bfrac{\log^2N}{\l N}.
\label{11q}\\
\nonumber\\ \nonumber\\
\ex[z_1^\prime-z_1\mid \bb,\bd]&=&\sum_{k\ge 2}\frac{kz_k}{2\mu}\,(k-1)\frac{2z_2}{2\mu}+
\e_{\ref{12}}.\label{12}\\
\ex[z_1^\prime-z_1\mid |\bv|]&=&\frac{z^2}{4\mu^2}\,\frac{\la^4f_0(\la)}{f_2(\la)^2}
+O\bfrac{\log^2N}{\l N}.
\label{12q}\\
\nonumber\\ \nonumber\\
\ex[y^\prime-y\mid \bb,\bd]&=&-1-\sum_{k\ge 3}\frac{ky_k}{2\mu}-
\sum_{k\ge 2}\frac{kz_k}{2\mu}\,(k-1)\frac{3y_3}{2\mu}+\e_{\ref{13}}.\label{13}\\
\ex[y^\prime-y\mid |\bv|]&=&-1-\frac{y}{2\mu}\,\frac{\la
 f_2(\la)}{f_3(\la)}-\frac{yz}{8\mu^2}\,\frac{\la^3}{f_3(\la)}\,
\,\frac{\la^2f_0(\la)}{f_2(\la)}+O\bfrac{\log^2N}{\l N}.
\label{13q}\\
\nonumber\\ \nonumber\\
\ex[z^\prime-z\mid \bb,\bd]&=&1-\sum_{k\ge 2}\frac{kz_k}{2\mu}-
\sum_{k\ge 2}\frac{kz_k}{2\mu}\,(k-1)\frac{2z_2}{2\mu}
+\sum_{k\ge 3}\frac{ky_k}{2\mu}+\e_{\ref{14}}.\label{14}\\
\ex[z^\prime-z\mid |\bv|]&=&1-\frac{z}{2\mu}\frac{\la
 f_1(\la)}{f_2(\la)}-\frac{z^2}{4\mu^2}\,\frac{\la^4f_0(\la)}{f_2(\la)^2}+\frac{y}{2\mu}
\,\frac{\la f_2(\la)}{f_3(\la)}+O\bfrac{\log^2N}{\l N}.
\label{14q}\\
\nonumber\\ \nonumber\\
\ex[\mu^\prime-\mu\mid \bb,\bd]&=&
-1-\sum_{k\ge 2}\frac{kz_k}{2\mu}(k-1)+\e_{\ref{15a}}.\label{15a}\\
\ex[\mu^\prime-\mu\mid |\bv|]&=&-1-\frac{z}{2\mu}\frac{\la^2 f_0(\la)}{f_2(\la)}
+O\bfrac{\log^2N}{\l N}.
\label{15aq}
\end{eqnarray}
\end{lemma}
\proofstart
The verification of \eqref{04x} -- \eqref{15a} is long but straightforward. We will verify
\eqref{04x} and
 \eqref{04xq} and add a few comments
 and hope that the reader is willing to accept or check the remainder by him/herself.

Suppose without loss of generality that \bx\ is such that $x_1=v=1\in Y_1$.
The remainder of \bx\ is a random
permutation of $2m-2\m$ $\star$'s and $2\m-1$ values from $[n]$ where the number
of times $j$ occurs is $\bd_\bx(j)$
for $j\in [n]$. The term -1 accounts for the deletion of $v$ from $\G$. There is a probability
$\frac{y_1}{2\m-1}=\frac{y_1}{2\m}+O\bfrac{1}{\m}$ that $x_2\in Y_1$ and this
accounts for the second term in \eqref{04x}.
Observe next that there is a probability $\frac{kz_k}{2\m-1}$ that $x_2\in Z_k,k\geq 2$.
In which case another $k-1$ edges will
be deleted. In expectation, the number of vertices in $Y_1$ lost by the deletion of
one such edge is $\frac{y_1-1}{2\m-3}$ and
this accounts for the third term. On the other hand, each such edge has a
$\frac{2y_2}{2\m-3}$ probability of being incident with
a vertex in $Y_2$. The deletion of such an edge will create a vertex in $Y_1$ and
this explains the fourth term. We collect the errors
from replacing $\m$ by $\m-1$ etc. into the last term. This gives a contribution of
order $1/N$. The above analysis ignored the extra contributions due to multiple edges.
We can bound this by
\beq{eta}
\eta_{\ref{04x}}=\sum_{k\geq 3}\frac{kz_k}{2\m}
\sum_{\ell\geq 3}\frac{\ell y_\ell}{2\m-1}\binom{k-1}{\ell-1}\bfrac{\ell}{2\m-k}^{\ell-2}.
\eeq
To explain this, we assume $x_2\in Z_k$, which is accounted for by the first sum over $k$. Now,
to create a vertex in $Y_1$, the removal of $x_2$ must
delete $\ell-1$ of the edges incident with some
vertex $y$ in $Y_\ell$. The term $\frac{\ell y_\ell}{2\m-1}$ is the probability that the first
of the chosen $\ell-1$ edges is incident with $y\in Y_\ell$ and the factor
$\bfrac{\ell}{2\m-k}^{\ell-2}$ bounds the probability that the remaining $\ell-2$ edges are
incident with $y$.

To go from conditioning on \bb,\bd\ to conditioning on $|\bv|$ we need to use the expected
values of $y_k,z_l$ etc., conditional on \bv. For this we use \eqref{f1} and \eqref{f2}.

We have, up to an error term $O\bfrac{\log^2N}{\l N}$,
\begin{eqnarray}
\ex\left[\sum_{k\ge 3}ky_k\biggr||\bv|\right]&=&\sum_{k\ge 3}k\,y\,\frac{\la^k}{k!f_3(\la)}
=\frac{y\la}{f_3(\la)}\sum_{j\ge 2}\frac{\la^j}{j!}=y\,\frac{\la f_2(\la)}{f_3(\la)},
\label{16}\\
\ex\left[\sum_{k\ge 2}kz_k\biggr||\bv|\right]&=&\sum_{k\ge 2}k\,z\,\frac{\la^k}{k!f_2(\la)}
=\frac{z\la}{f_2(\la)}\sum_{j\ge 1}\frac{\la^j}{j!}=z\,\frac{\la f_1(\la)}{f_2(\la)},\label{16a}\\
\ex\left[\sum_{k\ge 3}k(k-1)y_k\biggr||\bv|
\right]&=&\sum_{k\ge 3}k(k-1)\,y\,\frac{\la^k}{k!f_3(\la)}
=\frac{y\la^2}{f_3(\la)}\sum_{j\ge 1}\frac{\la^j}{j!}=y\,
\frac{\la^2f_1(\la)}{f_3(\la)},\label{17a}\\
\ex\left[\sum_{k\ge 2}k(k-1)z_k\biggr||\bv|
\right]&=&\sum_{k\ge 2}k(k-1)\,z\,\frac{\la^k}{k!f_2(\la)}
=\frac{z\la^2}{f_2(\la)}\sum_{j\ge 0}\frac{\la^j}{j!}=z\,\frac{\la^2f_0(\la)}{f_2(\la)}\label{17}.
\end{eqnarray}
In particular, using \eqref{17} in \eqref{04x} we get \eqref{04xq}. The other terms are obtained
in a similar fashion. We remark that we need to use \eqref{f2} when we deal with products
$z_ky_\ell$, $k\geq2$ and $\ell\ge3$.

Since, $k,\ell\leq\log n$ in \eqref{eta} we see, with the aid of \eqref{16} -- \eqref{17} that
$\ex[\eta_{\ref{04x}}\mid\bv]=O(1/N)$. This bound is true for all other $\e_?$.
\proofend

\subsection{Negative drift for $y_1,y_2,z_1$}
Algorithm \2G\ tries to keep $y_1,y_2,z_1$ small by its selection in Step 1. We now verify
that there is a negative drift in
$$\z=\z(t)=y_1+2y_2+z_1$$
in all cases of Step 1. This will enable us to show that \whp\ $\z$ remains small throughout the
execution of \2G. Let
\beq{eqQ}
Q=Q(\bv)=\frac{yz}{4\mu^2}\,\frac{\la^3}{f_3(\la)}\,\frac{\la^2f_0(\la)}{f_2(\la)}+
\frac{z^2}{4\mu^2}\,\frac{\la^4f_0(\la)}{f_2(\la)^2}.
\eeq
Then simple algebra gives
\begin{align}
&\ex[\z'-\z\mid|\bv|]=-(1-Q)-\brac{\z+y_2}\brac{\frac{1}{2\m}+\frac{z\l^2f_0(\l)}
{4\m^2f_2(\l)}}+O\bfrac{\log^2N}{\l N}&Case\ 1(a)\label{C1a}\\
&\ex[\z'-\z\mid|\bv|]=-2(1-Q)-\z\brac{\frac{1}{\m}+\frac{z\l^2f_0(\l)}
{2\m^2f_2(\l)}}+O\bfrac{\log^2N}{\l N}&Case\ 1(b)\label{C1b}\\
&\ex[\z'-\z\mid|\bv|]=-(1-Q)-\z\brac{\frac{1}{2\m}+\frac{z\l^2f_0(\l)}
{4\m^2f_2(\l)}}+O\bfrac{\log^2N}{\l N}&Case\ 1(c)\label{C1c}
\end{align}
We will show
\begin{lemma}\label{alem}[Pittel]
\beq{1-Q}
\l>0\text{ implies }Q<1
\eeq
and
\begin{equation}\label{Qless}
Q=\left\{\alignedat2
&O(\la^{-1}),\quad&&\la\to\infty,\\
&1-\Theta(\la^2),\quad&&\la\to 0.\endalignedat\right.
\end{equation}
\end{lemma}
\proofstart
Now, by \eqref{3}, $Q<1$ is equivalent to
$$
yz\frac{\la^5f_0(\la)}{f_2(\la)f_3(\la)}+z^2\frac{\la^4f_0(\la)}{f_2(\la)^2}
< \left(y\frac{\la f_2(\la)}
{f_3(\la)}+z\frac{\la f_1(\la)}{f_2(\la)}\right)^2,
$$
or, introducing $x=y/z$,
\begin{equation}\label{Fxla}
F(x,\la):=\frac{x\frac{\la^5f_0(\la)}{f_2(\la)f_3(\la)}+\frac{\la^4 f_0(\la)}
{f_2(\la)^2}}{\left(x\frac{\la f_2(\la)}{f_3(\la)}+
\frac{\la f_1(\la)}{f_2(\la)}\right)^2}<1,\quad\forall\,\l>0, x\ge 0.
\end{equation}
In particular, $F(\infty,\la)=0$. Now
$$
F_x(x,\la)=\left(x\frac{\la f_2(\la)}{f_3(\la)}+\frac{\la f_1(\la)}{f_2(\la)}\right)^{-4} G(x,\la),
$$
where
\begin{multline}\label{Gxla}
G(x,\la)=\frac{\la^5f_0(\la)}{f_2(\la)f_3(\la)}\left(x\frac{\la f_2(\la)}{f_3(\la)}+
\frac{\la f_1(\la)}{f_2(\la)}\right)^2\\
-2\left(x\frac{\la f_2(\la)}{f_3(\la)}+\frac{\la f_1(\la)}{f_2(\la)}\right)\frac{\la f_2(\la)}
{f_3(\la)}\left(x\frac{\la^5f_0(\la)}{f_2(\la)f_3(\la)}+
\frac{\la^4f_0(\la)}{f_2(\la)^2}\right).
\end{multline}
Notice that
$$
G(0,\la)=\la^6f_0(\la)f_1(\la)f_2(\la)^{-3}f_3(\la)^{-1}\bigl(\la f_1(\la)-2f_2(\la)\bigr)>0,
$$
as $\la f_1(\la)-2f_2(\la)>0$. Whence $F_x(0,\la)>0$ and
as a function of $x$,  $F(x,\la)$ attains its maximum at the root of $G(x,\la)=0$,
which is
\beq{bars}
\bar x=\frac{f_3(\la)\bigl(\la f_1(\la)-2f_2(\la)\bigr)}{\la f_2(\la)^2}.
\eeq
Now, \eqref{Gxla} implies that $\bar{x}$ satisfies
\beq{gf1}
\bar{x}\frac{\la^5f_0(\la)}{f_2(\la)f_3(\la)}+\frac{\la^4f_0(\la)}{f_2(\la)^2}=
\frac{\la^5f_0(\la)}{f_2(\la)f_3(\la)}
\left(\bar{x}\frac{\la f_2(\la)}{f_3(\la)}+\frac{\la f_1(\la)}{f_2(\la)}\right)
\times \frac{f_3(\la)}{2\la f_2(\la)}
\eeq
and \eqref{bars} implies that
\beq{gf2}
\bar{x}\frac{\la f_2(\la)}{f_3(\la)}+\frac{\la f_1(\la)}{f_2(\la)}=\frac{2(\l f_1(\l)-f_2(\l))}
{f_2(\l)}.
\eeq
Substituting \eqref{gf1}, \eqref{gf2} into \eqref{Fxla}, we see that
$$F(\bar x,\la)=\,\frac{\frac{\la^5f_0(\la)}{f_2(\la)f_3(\la)}\frac{f_3(\la)}{2\la f_2(\la)}}
{\left(\bar x\frac{\la f_2(\la)}{f_3(\la)}+\frac{\la f_1(\la)}{f_2(\la)}\right)}
=\,\frac{\la^4f_0(\la)}{4f_2(\la)\bigl(\la f_1(\la)-f_2(\la)\bigr)}.$$
Thus,
\begin{equation}\label{1-F}
1-F(\bar x,\la)=\frac{D(\la)}{4f_2(\la)\bigl(\la f_1(\la)-f_2(\la)\bigr)},
\end{equation}
where
\begin{align*}
D(\la)=&4f_2(\la)\bigl(\la f_1(\la)- f_2(\la)\bigr)-\la^4f_0(\la)\\
=&-4-4\la-(\la^4+4\la^2-8)e^{\la}+(4\la-4)e^{2\la}.
\end{align*}
In particular,
\begin{equation}\label{1-F,inf}
1-F(\bar x,\la)=1-O(\la^{-1}),\quad\la\to\infty.
\end{equation}
Expanding $e^{\la}$and $e^{2\la}$, we obtain after collecting like terms that
$$D(\la)=\sum_{j\ge 6}\frac{d_j}{j!}\,\la^j,$$
where
$$d_j=2^{j+1}(j-2)-(j)_4-4(j)_2+8.$$
Here $d_j=0$ for $0\leq j\leq 5$ and $d_6=40,d_7=280,d_8=1176,d_9=3864,d_{10}=10992$ and $d_j>0$ for $j\ge 11$ is clear.
Therefore $D(\la)$ is positive for all $\la>0$. Since $D(\la)\sim d_6\la^6$
and $4f_2(\la)\bigl(\la f_1(\la)-f_2(\la)\sim \l^4$
as $\la\to 0$,  we see that
\begin{equation}\label{1-F,zero}
1-F(\bar x,\la)\sim d_6\la^2,\quad \la\to 0.
\end{equation}
This completes the proof of Lemma \ref{alem}.
\proofend

It follows from \eqref{C1a}, \eqref{C1b}, \eqref{C1c} and Lemma \ref{alem} that,
regardless of case,
\beq{eqx1}
\z>0\text{ implies }\ex[\z'-\z\mid|\bv|]\leq -c_1(1\wedge\l)^2+O\bfrac{\log^2N}{\l N}
\eeq
for some absolute constant $c_1>0$, where $\wel=\min\set{1,\l}$.

To avoid dealing with the error term in \eqref{eqx1} we introduce the stopping time,
$$T_{er}=\min\set{t:\l^2\leq\frac{\log^3n}{\l N}}.$$
(This is well defined, since eventually $N=0$).

The following stopping time is also used:
$$T_0=\min\set{t:\l\leq1\text{ or }
N\leq n/2}<T_{er}.$$

So we can replace \eqref{eqx1} by
\beq{eqx2}
\z>0\text{ implies }\ex[\z'-\z\mid|\bv|]\leq -c_1/2,\qquad 0\leq t\leq T_0,
\eeq
which holds for $n$ sufficiently large.

There are several places where we need a bound on $\l$:
\begin{lemma}\label{lambda}
{\bf Whp} $\l\leq 3ce$ for $t\leq T_0$.
\end{lemma}
\proofstart
We will show that w.h.p. $y_1+2y_2+z_1=o(n)$ throughout.
It follows from \eqref{3} and the inequalities in Section \ref{simpineq} that if
$\La$ is sufficiently large and if $\l(t)\geq \La$ then $Y\cup Z$
contains $y+z$ vertices and at least $\La(y+z)/2$ edges and hence has total
degree at least $\La(y+z)$.
We argue that \whp\ $G$ does not contain such a sub-graph.
We will work in the random sequence model. We can assume that $|Y\cup Z|\geq n/2$.
Now fix a set $S\subseteq [n]$ where $s=|S|\geq n/3$. Let
$D$ denote the total degree of vertices in $S$. Then
\begin{multline}\label{boundd}
\Pr(D=d)\leq O(n^{1/2})\sum_{\substack{d_1+\cdots+d_s=d\\d_j\geq 3}}
\prod_{i=1}^s\frac{\l^{d_i}}{f_3(\l)d_i!}\leq O(n^{1/2})\frac{\l^d}{d!f_3(\l)^s}
\sum_{\substack{d_1+\cdots+d_s=d\\d_j\geq0}}\frac{d!}{d_1!\cdots d_s!}\\
=O(n^{1/2})\frac{\l^ds^d}{d!f_3(\l)^s}.
\end{multline}
Here $\l=\l(0)$ and we are using Lemma \ref{lem3}. The factor $O(n^{1/2})$ accounts for the
conditioning that the total degree is $2cn$. Now $\l(0)\leq 2c$ and $f_3(\l(0))\geq 1$. It follows
that
$$
\Pr(\exists S:d\geq \La s)\leq O(n^{1/2})
\sum_{s\geq n/3}\sum_{d\geq \La s}\binom{n}{s}
\frac{(2c)^ds^d}{d!}\leq
O(n^{1/2})\sum_{s\geq n/3}\sum_{d\geq \La s}\bfrac{ne}{s}^s
\frac{(2c)^ds^d}{d!}
$$
The terms involving $d$ in the second sum are $u_d=\frac{(2cs)^d}{d!}$ and for $d/s$ large we
have $u_{d+1}/u_d=O(s/d)$ and so we can put $d=\La s$ in the second expression.
After substituting $d!\geq (d/e)^d$ this gives
$$\Pr(\exists S:d\geq \La s)\leq O(n^{1/2})\sum_{s\geq n/2}
\bfrac{3e(2ce)^{\La}}{\La^{\La}}^s=o(1)
$$
if $\La\geq 3ce$.
\proofend

Our aim now is to give a high probability bound on the maximum value that $\z$ will take
during the process.
We first prove a simple lemma involving the functions $\f_j(x)=\frac{xf_{j-1}(x)}{f_j(x)}$, $j=2,3$.
\begin{lemma}\label{fuk}
\beq{A1}
\f_j(x)\text{ is convex and increasing and }
j\leq \f_j(x)\text{ and }\frac{1}{j+1}\leq \f_j'(x)\leq 1\text{ for }j=2,3.
\eeq
\end{lemma}
\proofstart
Now, if $H(x)=\frac{xF(x)}{G(x)}$ then
$$H'(x)=\frac{G(x)(xF'(x)+F(x))-xF(x)G'(x)}{G(x)^2}$$
and
\begin{multline*}
H''(x)=\\
\frac{2xF(x)G'(x)^2+G(x)^2(2F'(x)+xF''(x))-G(x)(2xF'(x)G'(x)+F(x)(2G'(x)+xG''(x)))}{G(x)^3}.
\end{multline*}
{\bf Case $j=2$:}
\beq{deeriv}
\f_2'(x)=\frac{e^{2x}-(x^2+2)e^x+1}{(e^x-1-x)^2}.
\eeq
But,
$$e^{2x}-(x^2+2)e^x+1=\sum_{j\geq 4}\frac{2^j-j(j-1)-2}{j!}x^j$$
and so $\f_2'(x)>0$ for $x>0$.
\beq{deeeriv}
\f_2''(x)=\frac{e^{2x}(x^2-4x+2)+e^x(x^3+x^2+4x-4)+2}{(e^x-1-x)^3}.
\eeq
But
\begin{multline*}
e^{2x}(x^2-4x+2)+e^x(x^3+x^2+4x-4)+2=\\
\sum_{j\geq 6}\frac{2^{j-2}(j(j-1)-8j+8)+j(j-1)(j-2)+j(j-1)+4j-4}{j!}x^j
\end{multline*}
and so $\f_2''(x)>0$ for $x>0$.

{\bf Case $j=3$:}
\beq{deeriv0}
\f_3'(x)=\frac{2e^{2x}-e^x\brac{x^3-x^2+4x+4}+x^2+4x+2}{2(e^x-1-x-\frac{x^2}{2})^2}.
\eeq
But,
$$2e^{2x}-e^x\brac{x^3-x^2+4x+4}+x^2+4x+2=\sum_{j\geq 6}\frac{2^{j+1}-j(j-1)(j-2)+j(j-1)-4j-4}{j!}x^j$$
and so $\f_3'(x)>0$ for $x>0$.
\beq{deeeriv1}
\f_3''(x)=\frac{x(e^{2x}(2x^2-12x+12)+e^x(x^4+8x^2-24)+2x^2+12x+12)}{4(e^x-1-x-\frac{x^2}{2})^3}.
\eeq
But
\begin{multline*}
e^{2x}(2x^2-12x+12)+e^x(x^4+8x^2-24)+2x^2+12x+12\\
\sum_{j\geq 9}\frac{2^{j-1}(j(j-1)-12j+24)+j(j-1)(j-2)(j-3)+8j(j-1)-24}{j!}x^j.
\end{multline*}
and so $\f_3''(x)>0$ for $x>0$.

So $\f_2,\f_3$ are convex and so we only need to check that
$\f_2(0)=2,\f_2'(0)=1/3,\f_3(0)=3,\f_3'(0)=1/4$ and $\f_2'(\infty)=\f_3'(\infty)=1$.
\proofend

Consider $\l$ as a function of $\bv$, defined by
\beq{2x}
y\f_3(\l)+z\f_2(\l)=\Pi
\eeq
where $\Pi=2\m-y_1-2y_2-z_1$.

We now prove a
lemma bounding the change in $\l$ as we change \bv.
\begin{lemma}\label{onestep}
$$|\l(\bv_1)-\l(\bv_2)|= O\bfrac{||\bv_1-\bv_2||_1}{N},\qquad\text{ for }t<T_{er}.$$
\end{lemma}
\proofstart
We write $\bv_1=(y_1,y_2,z_1,y,z,\m)\geq 0$ and $\bv_2=(y_1+\d_{y_1},y_2+\d_{y_2},z_1+\d_{z_1},
y+\d_y,z+\d_z,\m+d_{\m})\geq 0$ and $\Pi,\Pi+\d_\Pi$ for the two values of $\Pi$.
Then
\beq{3x}
(y+\d_y)\f_3(\l+\d_\l)-y\f_3(\l)+(z+\d_ z)\f_2(\l+\d_\l)-\f_2(\l)=\d_\Pi.
\eeq
Convexity and our lower bound on $\f_j'$ implies that
$$\f_j(\l)\geq \f_j(\l+\d_\l)-\d_\l\f_j'(\l+\d_\l)\geq \f_j(\l+\d_\l)-\d_\l.$$
So from \eqref{3x} we have
$$(y+\d_y)(\f_3(\l)+\d_\l)-y\f_3(\l)+(z+\d_z)(\f_2(\l)+\d_\l)-\f_2(\l)\geq\d_\Pi.$$
This implies that
$$\d_\l\geq \frac{\d_\Pi-\d_y\f_3(\l)-\d_z\f_2(\l)}{y+\d_y+z+\d_z}.$$
So,
$$\d_\l\leq 0\text{ implies }|\d_\l|=O\bfrac{||\bv_1-\bv_2||_1}{N}.$$
Note that we use Lemma \ref{lambda} to argue that $\f_j(\l),j=2,3$ are bounded
within our range of interest.

To deal with $\d_\l\geq 0$ we observe that convexity implies
$$\f_j(\l+\d_\l)\geq \f_j(\l)+\d_\l\f_j'(\l).$$
So from \eqref{3x} we have
$$(y+\d_y)(\f_3(\l)+\d_\l\f_3'(\l))-y\f_3(\l)+(z+\d_z)(\f_2(\l)+\d_\l\f_2'(\l))-\f_2(\l)\leq\d_\Pi.$$
This implies that
$$\d_\l\leq \frac{\d_\Pi-\d_y\f_3(\l)-\d_z\f_2(\l)}{(y+\d_y)\f_3'(\l)+(z+\d_z)\f_2'(\l)}.$$
So,
$$\d_\l\geq 0\text{ implies }|\d_\l|= O\bfrac{||\bv_1-\bv_2||_1}{N}.$$
\proofend

\begin{lemma}\label{Klem}
If $c\geq 15$ then \qs
$$\not\exists 1\leq t\leq T_0:\;\z(t)> \log^2 n.$$
\end{lemma}
\proofstart
Define a sequence
$$X_i=\begin{cases}\min\set{\z(i+1)-\z(i),\log n}
&0\leq i\leq T_0\\-c_1/2&T_0<i\leq n\end{cases}$$
The variables $X_1,X_2,\ldots,X_n$ are not independent. On the other hand, conditional on an
event that occurs \qs, we see that
$$X_{s+1}+\ldots+X_t=\z(t)-\z(s)\text{ for }0\leq s<t\leq T_0$$
and
$$\ex[X_t\mid X_1,\ldots,X_{t-1}]\leq -c_1/2\text{ for }t\leq n.$$
Next, for $0\leq s\leq t\leq T_0$ let
$$\lb(s,t)=\sum_{\t=s+1}^t\l(\t)^2.$$
Note that
\beq{ign}
\lb(s,t)\geq t-s.
\eeq
We argue as in the proof of the Azuma-Hoeffding inequality that for any $1\leq s<t\leq n$
and $u\geq 0$,
\beq{x4}
\Pr(X_{s+1}+\cdots+X_t\geq u-c_1\lb(s,t)/2)\leq \exp\set{-\frac{2u^2}{(t-s)\log^2n}}.
\eeq
We deduce from this that
\begin{multline}\label{eqx4}
\Pr(\exists 1\leq s<t\leq T_0:\z(s)=0<\z(\t),s<\t\leq t)\leq\\
n^2\exp\set{-\frac{2\max\set{0,c_1\lb(s,t)/2-\log n}^2}{(t-s)\log^2n}}.
\end{multline}
Putting
$t-s=L_1=\log^2n$ we see from \eqref{eqx4} that \qs
\beq{eqx5}
\not\exists 1\leq s<t-L_1\leq T_0-L_1:\z(s)=0<\z(\t),s<\t\leq t.
\eeq
Suppose now that there exists $\t\leq T_0$ such that $\z(\t)\geq L_1$. Then q.s. there exists
$t_1\leq \t\leq t_1+L_1$ such that $\z(t_1)=0$. But then given $t_1$,
$$\Pr(\exists t_1\leq \t\leq t_1+L_1:\z(\t)\geq L_1)\leq \exp\set{-\frac{2(c_1L_1/2-\log n)^2}{L_1\log^2n}}.$$
Here we are using the generalisation of Hoeffding-Azuma that deals with $\max_{i\leq L_1}X_1+\cdots+X_i$.

And then we get that \qs
\beq{eqx6}
\not\exists t\leq T_0:\z(\t)\geq L_1.
\eeq
We do this in two stages because of the condition $\z>0$ in \eqref{eqx2}. Remember here
that $\z(0)=0$ and \eqref{eqx5} says that $\z$ cannot stay positive for very long.
\proofend

\section{Associated Equations.}\label{diff}
The expected changes conditional on \bv\ lead us to consider the following collection of
differential equations: Note that we do not use any scaling. We will put hats on variables
i.e. $\hy_1$ etc. will be the deterministic counterpart of $y_1$. Also,
as expected, the hatted equivalent of \eqref{2x} holds:
\beq{10xxx}
\frac{\hy\hl f_2(\hl)}{f_3(\hl)}+\frac{\hz\l \f_1(\hl)}{f_2(\hl)}=2\hm-\hy_1-2\hy_2-\hz_1.
\eeq

\noindent
{\bf Step 1(a).\/} $\hy_1>0$.
\begin{align}
\frac{d\hy_1}{dt}=&\,-1-\frac{\hy_1}{2\hm}-\frac{\hy_1\hz}{4\hm^2}\frac{\hl^2f_0(\hl)}{f_2(\hl)}
+\frac{\hy_2\hz}{2\hm^2}\,\frac{\hl^2f_0(\hl)}{f_2(\hl)},\label{18x}\\
\frac{d\hy_2}{dt}=&\,-\frac{\hy_2}{\hm}-\frac{\hy_2\hz}{2\hm^2}\frac{\hl^2f_0(\hl)}{f_2(\hl)}
+\frac{\hy\hz}{8\hm^2}\frac{\hl^3}{f_3(\hl)}\,\frac{\hl^2f_0(\hl)}{f_2(\hl)},\label{18}\\
\frac{d\hz_1}{dt}=&\,-\frac{\hz_1}{2\hm}-\frac{\hz_1\hz}{4\hm^2}\,\frac{\hl^2f_0(\hl)}{f_2(\hl)}
+\frac{\hz^2}{4\hm^2}\frac{\hl^4f_0(\hl)}{f_2(\hl)^2},\label{19}\\
\frac{d\hy}{dt}=&\,-\frac{\hy}{2\hm}\,\frac{\hl
 f_2(\hl)}{f_3(\hl)}-\frac{\hy\hz}{8\hm^2}\,\frac{\hl^3}{f_3(\hl)}\,
\,\frac{\hl^2f_0(\hl)}{f_2(\hl)},\label{20}\\
\frac{dz}{dt}=&\,\frac{\hy}{2\hm}\frac{\hl f_2(\hl)}{f_3(\hl)}-\frac{\hz}{2\hm}\frac{\hl
 f_1(\hl)}{f_2(\hl)}-\frac{\hz^2}{4\hm^2}\,\frac{\hl^4f_0(\hl)}{f_2(\hl)^2},\label{21}\\
\frac{d\hm}{dt}=&\,-1-\frac{\hz}{2\hm}\frac{\hl^2f_0(\hl)}{f_2(\hl)}.\label{22}
\end{align}
{\bf Step 1(b).\/} $\hy_1=0,\hy_2>0$.
\begin{align}
\frac{d\hy_1}{dt}=&\,\frac{\hy_2\hz}{\hm^2}\,\frac{\hl^2f_0(\hl)}{f_2(\hl)},\label{18xx}\\
\frac{d\hy_2}{dt}=&\,-1-\frac{2\hy_2}{\hm}-\frac{\hy_2\hz}{\hm^2}\frac{\hl^2f_0(\hl)}{f_2(\hl)}
+\frac{\hy\hz}{4\hm^2}\frac{\hl^3}{f_3(\hl)}\,\frac{\hl^2f_0(\hl)}{f_2(\hl)},\label{18mm}\\
\frac{d\hz_1}{dt}=&\,-\frac{\hz_1}{\hm}-\frac{\hz_1\hz}{2\hm^2}\,\frac{\hl^2f_0(\hl)}{f_2(\hl)}
+\frac{\hz^2}{2\hm^2}\frac{\hl^4f_0(\hl)}{f_2(\hl)^2},\label{19mm}\\
\frac{d\hy}{dt}=&\,-\frac{\hy}{\hm}\,\frac{\hl
 f_2(\hl)}{f_3(\hl)}-\frac{\hy\hz}{4\hm^2}\,\frac{\hl^3}{f_3(\hl)}\,
\,\frac{\hl^2f_0(\hl)}{f_2(\hl)},\label{20mm}\\
\frac{d\hz}{dt}=&\,\frac{\hy}{\hm}\frac{\hl f_2(\hl)}{f_3(\hl)}-\frac{\hz}{\hm}\frac{\hl
 f_1(\hl)}{f_2(\hl)}-\frac{\hz^2}{2\hm^2}\,\frac{\hl^4f_0(\hl)}{f_2(\hl)^2},\label{21mm}\\
\frac{d\hm}{dt}=&\,-2-\frac{\hz}{\hm}\frac{\hl^2f_0(\hl)}{f_2(\hl)}.\label{22mm}
\end{align}

{\bf Step 1(c).\/} $\hy_1=\hy_2=0,\hz_1>0$.

\begin{align}
{ \frac{d\hy_1}{dt}=}&0,\label{23x}\\
{ \frac{d\hy_2}{dt}=}&\,{\frac{\hy\hz}{8\hm^2}\,\frac{\hl^3}{f_3(\hl)}\,
\frac{\hl^2f_0(\hl)}{f_2(\hl)}},\label{23}\\
{ \frac{d\hz_1}{dt}=}&{ \,-1-\frac{\hz_1}{2\hm}-\frac{\hz_1\hz}{4\hm^2}\,\frac{\hl^2f_0(\hl)}{f_2(\hl)}+
\frac{\hz^2}{4\hm^2}\,\frac{\hl^4f_0(\hl)}{f_2(\hl)^2}},
\label{24}\\
{ \frac{d\hy}{dt}=}&{\,-\frac{\hy}{2\hm}\,
\frac{\hl f_2(\hl)}{f_3(\hl)}-\frac{\hy\hz}{8\hm^2}\,\frac{\hl^3}{f_3(\hl)}\,
\,\frac{\hl^2f_0(\hl)}{f_2(\hl)},}\label{20a}\\
\frac{d\hz}{dt}=&\,\frac{\hy}{2\hm}\,\frac{\hl f_2(\hl)}{f_3(\hl)}
-\frac{\hz}{2\hm}\frac{\hl f_1(\hl)}{f_2(\hl)}
-\frac{\hz^2}{4\hm^2}\,\frac{\hl^4f_0(\hl)}{f_2(\hl)^2},\label{21a}\\
\frac{d\hm}{dt}=&\,-1-\frac{\hz}{2\hm}\frac{\hl^2f_0(\hl)}{f_2(\hl)}.\label{25a}
\end{align}

{\bf Step 2.\/} $\hy_1=\hy_2=\hz_1=0$.
\begin{align}
\frac{d\hy_1}{dt}=&0,\label{18ax}\\
\frac{d\hy_2}{dt}=&\,\frac{\hy\hz}{8\hm^2}
\frac{\hl^3}{f_3(\hl)}\frac{\hl^2f_0(\hl)}{f_2(\hl)},\label{18a}\\
\frac{d\hz_1}{dt}=&\,\frac{\hz^2}{4\hm^2}\,\frac{\hl^4f_0(\hl)}{f_2(\hl)^2},\label{19a}\\
\frac{d\hy}{dt}=&\,-1-\frac{\hy}{2\hm}\,\frac{\hl
 f_2(\hl)}{f_3(\hl)}-\frac{\hy\hz}{8\hm^2}\,\frac{\hl^3}{f_3(\hl)}\,
\,\frac{\hl^2f_0(\hl)}{f_2(\hl)},\label{20aa}\\
\frac{d\hz}{dt}=&\,1-\frac{\hz}{2\hm}\frac{\hl
 f_1(\hl)}{f_2(\hl)}-\frac{\hz^2}{4\hm^2}\,\frac{\hl^4f_0(\hl)}{f_2(\hl)^2}+\frac{\hy}{2\hm}
\,\frac{\hl f_2(\hl)}{f_3(\hl)},\label{21aa}\\
\frac{d\hm}{dt}=&\,-1-\frac{\hz}{2\hm}\frac{\hl^2 f_0(\hl)}{f_2(\hl)}.\label{22aa}
\end{align}
We will show that \whp\ the process defined by \2G\ can be closely modeled by a suitable
weighted sum of the above four sets of equations.
Let these weights be $\th_a,\th_b,\th_c$ and $1-\th_a-\th_b-\th_c$ respectively.
It has been determined that $y_1,y_2,z_1$ are
all $O(\log^2 n)$ \whp. We will only need to analyse our process up till the time $y=0$ and
we will show that at this time, $z=\Omega(n)$ \whp. Thus $y_1,y_2,z_1$ are "negligible" throughout.
In which case $\hy_1,\hy_2,\hz_2$ should also be negligible.
It makes sense therefore to choose $\th_a=0$.
The remaining weights should be chosen so that the
weighted derivatives of $\hy_1,\hy_2,\hz_1$ are zero.
This has all been somewhat heuristic and its validity
will be verified in Section \ref{section-close}.
\subsection{Sliding trajectory}
Conjecturally we need to mix Steps 1(a), 1(b) 1(c) and 2 with
nonnegative weights $\theta_a=0$,
$\theta_b$, $\theta_c$, $\th_2=1-\theta_b-\theta_c$ respectively, chosen such that
the resulting system of differential equations admits a solution
such that $\hy_2(t)\equiv 0$ and $\hz_1(t)
\equiv 0$.

We will write the multipliers in terms of
\beq{ABCD}
\hA=\frac{\hy\hz\hl^5f_0(\hl)}{8\hm^2f_2(\hl)f_3(\hl)},\quad \hB=\frac{\hz^2\hl^4f_0(\hl)}
{4\hm^2f_2(\hl)^2},
\quad \hC=\frac{\hy\hl f_2(\hl)}{2\hm f_3(\hl)},\quad \hD=\frac{\hz\hl^2f_0(\hl)}{2\hm f_2(\hl)}.
\eeq
Using, \eqref{18x}, \eqref{18xx}, \eqref{23x} and \eqref{18ax} we see $\hy_1(t)\equiv 0$ implies that
$$0\equiv\frac{d\hy_1}{dt}=\th_a.$$
Equivalently
\begin{equation}\label{thetaaa}
\theta_a=0.
\end{equation}

Using \eqref{18},
\eqref{23} and \eqref{18a}, we see that $\hy_2(t)\equiv 0$ implies that
\begin{align}
&0\equiv\frac{d\hy_2}{dt}\nonumber\\
&=\theta_b\left[-1+\frac{\hy\hz}{4\hm^2}\,\frac{\hl^3}{f_3(\hl)}\,\frac{\hl^2f_0(\hl)}{f_2(\hl)}\right]
+\theta_c\frac{\hy\hz}{8\hm^2}\,\frac{\hl^3}{f_3(\hl)}\,\frac{\hl^2f_0(\hl)}{f_2(\hl)}
+(1-\theta_b-\theta_c)\frac{\hy\hz}{8\hm^2}\,\frac{\hl^3}{f_3(\hl)}\,\frac{\hl^2f_0(\hl)}{f_2(\hl)},\nonumber\\
&=-(1-\hA)\th_b+\hA.\label{dy2=0}
\end{align}
Equivalently
\begin{equation}\label{thetaa}
\theta_b=\frac{\hA}
{1-\hA}.
\end{equation}
Likewise, using \eqref{19}, \eqref{24} and \eqref{19a}, $z_1(t)\equiv 0$ implies
\begin{align}
&0\equiv\frac{d\hz_1}{dt}\nonumber\\
&=\theta_b \frac{\hz^2}{2\hm^2}\,\frac{\hl^4f_0(\hl)}{f_2(\hl)^2}
+\theta_c\left[-1+\frac{\hz^2}{4\hm^2}\,\frac{\hl^4f_0(\hl)}{f_2(\hl)^2}\right]
+(1-\theta_b-\theta_c)\frac{\hz^2}{4\hm^2}\,\frac{\hl^4f_0(\hl)}{f_2(\hl)^2},\nonumber\\
&=\hB\th_b-\th_c+\hB.\label{dz1=0}
\end{align}
Equivalently
\begin{equation}\label{thetab}
\theta_c=(1+\theta_b)\hB=\frac{\hB}{1-\hA}.
\end{equation}
From \eqref{thetaa} it follows that $\theta_b\ge 0$ iff
$$
\hA\leq  1,
$$
in which case, by \eqref{thetab}, $\theta_c\ge 0$, as well.
From \eqref{dz1=0} and \eqref{thetab} it follows that
$1-\theta_b-\theta_c\ge 0$ iff
\beq{2A}
2\hA+\hB\leq 1.
\eeq
We conclude that $\theta_b,\theta_c,1-\theta_b-\theta_c\in [0,1]$ iff $Q\leq 1$, see \eqref{eqQ}.
But this is implied by Lemma \ref{alem}.

It may be of some use to picture the equations defining $\th_a,\th_b,\th_c,\th_2$:
\beq{equations}
\begin{array}{rrrrc}
-\th_a&&&&=0\\
&(1-\hA)\th_b&&&=\hA\\
&-\hB\th_b&+\th_c&&=\hB\\
\th_a&+\th_b&+\th_c&+\th_2&=1.
\end{array}
\eeq

If in the notation of Lemma \ref{Klem} we
let $\Omega_1=\set{\bv:\z\leq L_1}$ then we may restrict our attention to \bv\
in \eqref{04x} -- \eqref{15a} such that $\bv\in\Omega_1$.
In which case, the terms involving $y_1,y_2,z$ can be absorbed
into the error term fopr $t\leq T_0$. The relevant equations then become, with
$$A=\frac{yz\l^5f_0(\l)}{8\m^2f_2(\l)f_3(\l)},\quad B=\frac{z^2\l^4f_0(\l)}
{4\m^2f_2(\l)^2},
\quad C=\frac{y\l f_2(\l)}{2\m f_3(\l)},\quad D=\frac{z\l^2f_0(\l)}{2\m f_2(\l)}.$$
{\bf Step 1(a)\/}. $y_1>0$.
\begin{eqnarray}
\ex[y_1^\prime-y_1\mid
|\bv|]&=&-1+O\bfrac{\log^2N}{\l N}\label{04xq1mmm}\\
\ex[y_2^\prime-y_2\mid |\bv|]&=&A+O\bfrac{\log^2N}{\l N}\label{04q1mmm}\\
\ex[z_1^\prime-z_1\mid |\bv|]&=&B+O\bfrac{\log^2N}{\l N}.\label{05q1mm}\\
\ex[y^\prime-y\mid |\bv|]&=&- C-A+O\bfrac{\log^2N}{\l N}.\label{06q1}\\
\ex[z^\prime-z\mid |\bv|]&=&C-(1-C)-B +O\bfrac{\log^2N}{\l N}.\label{07q1}\\
\ex[\mu^\prime-\mu\mid|\bv|]&=&-1-D+O\bfrac{\log^2N}{\l N}.\label{08q1}
\end{eqnarray}
{\bf Step 1(b)\/}. $y_1=0,y_2>0$.
\begin{eqnarray}
\ex[y_1^\prime-y_1\mid
 |\bv|]&=&O\bfrac{\log^2N}{\l N}.\label{4xq1mm}\\
\ex[y_2^\prime-y_2\mid
 |\bv|]&=&-1+2A+O\bfrac{\log^2N}{\l N}.\label{4q1mm}\\
\ex[z_1^\prime-z_1\mid
|\bv|]&=&2B+O\bfrac{\log^2N}{\l N}.\label{5q1}\\
\ex[y^\prime-y\mid |\bv|]&=&-2C-2A+O\bfrac{\log^2N}{\l N}.\label{6q1}\\
\ex[z^\prime-z\mid |\bv|]&=&2C-2(1-C)-2B+O\bfrac{\log^2N}{\l N}.\label{7q1}\\
\ex[\mu^\prime-\mu\mid |\bv|]&=&-2-2D+O\bfrac{\log^2N}{\l N}.\label{8aq1}
\end{eqnarray}
{\bf  Step 1(c).\/} $y_1=y_2=0,z_1>0$.
\begin{eqnarray}
\ex[y_1^\prime-y_1\mid |\bv|]&=&O\bfrac{1}{N}.\label{9xq1mmm}\\
\ex[y_2^\prime-y_2\mid |\bv|]&=&A+O\bfrac{\log^2N}{\l N}.\label{9q1mmm}\\
\ex[z_1^\prime-z_1\mid |\bv|]&=&-1+B+O\bfrac{\log^2N}{\l N}.\label{10q1mmm}\\
\ex[y^\prime-y\mid |\bv|]&=&-C-A+O\bfrac{\log^2N}{\l N}.\label{6aq1}\\
\ex[z^\prime-z\mid |\bv|]&=&C-(1-C)-B+O\bfrac{\log^2N}{\l N}.\label{7aq1}\\
\ex[\mu^\prime-\mu\mid |\bv|]&=&-1-D+O\bfrac{\log^2N}{\l N}.\label{8aaq1}
\end{eqnarray}

{\bf Step 2.\/} $y_1=y_2=z_1=0$.
\begin{eqnarray}
\ex[y_1^\prime-y_1\mid |\bv|]&=&O\bfrac{1}{N}.\label{11xq1}\\
\ex[y_2^\prime-y_2\mid |\bv|]&=&A+O\bfrac{\log^2N}{\l N}.
\label{11q1}\\
\ex[z_1^\prime-z_1\mid |\bv|]&=&B+O\bfrac{\log^2N}{\l N}.
\label{12q1}\\
\ex[y^\prime-y\mid |\bv|]&=&-1-C-A+O\bfrac{\log^2N}{\l N}.
\label{13q1}\\
\ex[z^\prime-z\mid |\bv|]&=&1+C-(1-C)-B+O\bfrac{\log^2N}{\l N}.
\label{14q1}\\
\ex[\mu^\prime-\mu\mid |\bv|]&=&-1-D+O\bfrac{\log^2N}{\l N}.
\label{15aq1}
\end{eqnarray}
\subsection{Closeness of the process and the differential equations}\label{section-close}
We already know that $y_1,y_2,z_1$ are small \whp\ up to
time $T_0$. We now show that \whp\ $y,z,\mu$ are closely
approximated by $\hy,\hz,\hm$, which are the solutions to
the weighted sum of the sets of equations labelled
Step 1(b), Step 1(c) and Step 2. These equations will be
simplified by putting $y_1=y_2=z_1=0$. First some notation.
We will use $\psi_{\xi,\eta}$ to denote the expression we
have obtained for the derivative of $\xi$ in Case 1 ($\eta$)
or Case 2 in the case of $\eta=2$.
We are then led to
consider the equations:

{\bf Sliding Trajectory:}
\begin{align*}
\frac{d\hy}{dt}&=\th_b\psi_{b,y}(\hy,\hz,\hm)+
\th_c\psi_{c,y}(\hy,\hz,\hm)+(1-\th_b-\th_c)\psi_{2,y}(\hy,\hz,\hm)\\
&=\th_b(-2(\hC+\hA))+\th_c(-(\hC+\hA))+(1-\th_b-\th_c)(-(1+\hC+\hA))\\
&=-(\hC+\hA)(2\th_b+\th_c+1-\th_b-\th_c)-(1-\th_b-\th_c)\\
&=\frac{\hB-\hC}{1-\hA}-1.\\
\frac{d\hz}{dt}&=\th_b\psi_{b,z}(\hy,\hz,\hm)+
\th_c\psi_{c,z}(\hy,\hz,\hm)+(1-\th_b-\th_c)\psi_{2,z}(\hy,\hz,\hm)\\
&=\th_b(2(\hC-(1-\hC)-\hB))+\th_c(\hC-(1-\hC)-\hB)+(1-\th_b-\th_c)(1+\hC-1+\hC-\hB)\\
&=(2\hC-\hB)(\th_b+1)-2\th_b-\th_c\\
&=\frac{2\hC-2\hA-2\hB}{1-\hA}.\\
\frac{d\hm}{dt}&=\th_b\psi_{b,\m}(\hy,\hz,\hm)+
\th_c\psi_{c,\m}(\hy,\hz,\hm)+(1-\th_b-\th_c)\psi_{2,\m}(\hy,\hz,\hm)\\
&=\th_b(-2(1+\hD))+\th_c(-(1+\hD))+(1-\th_b-\th_c)(-(1+\hD))\\
&=-(1+\hD)(2\th_b+\th_c+1-\th_b-\th_c)\\
&=-\frac{1+\hD}{1-\hA}.
\end{align*}
The starting conditions are
\beq{hinitial}
\hy(0)=n,\hz(0)=0,\hm(0)=cn.
\eeq
Summarising:
\beq{slide}
\frac{d\hy}{dt}=\frac{\hB-\hC}{1-\hA}-1;\quad\frac{d\hz}{dt}=\frac{2\hC-2\hA-2\hB}{1-\hA};\quad\frac{d\hm}{dt}=-\frac{1+\hD}{1-\hA}.
\eeq
and
\beq{slide1}
\frac{\hy\hl f_2(\hl)}{f_3(\hl)}+\frac{\hz\hl f_1(\hl)}{f_2(\hl)}=2\hm.
\eeq
We remark for future reference that \eqref{slide} implies that
\beq{future}
\hm\text{ is decreasing with $t$ as long as $\hl>0$}
\eeq
and \eqref{slide1} implies that
\beq{slide2}
\hy+\hz\leq \frac{2\hm}{\hl}.
\eeq

Let $\bu=\bu(t)$ denote $(y(t),z(t),\m(t))$ and let $\hu=\hu(t)$ denote $(\hy(t),\hz(t),\hm(t))$.
We now show that $\bu$ and $\hu$ remain close:
\begin{lemma}\label{close}
$$||\bu(t)-\hu(t)||_1\leq  n^{8/9}, \qquad\text{ for }1\leq t\leq T_0,\ \whp.$$
\end{lemma}
\proofstart
Let $\d_\eta(\bv),\eta=a,b,c,2$ be the 0/1 indicator for the process
\2G\ applying
Step 1($\eta$) for $\eta=a,b,c$ or Step 2 if $\eta=2$ when the current state is \bv.
For times $t_1<t_2$ we use the notation
$$\D_\eta(\bv(t_1,t_2))=\sum_{t=t_1}^{t_2}\d_\eta(\bv(t))$$

Now let $\r=n^\a$ where $\a=1/4$. It follows from Lemma \ref{onestep} that for $t\leq T_0-\r$,
\beq{pp0}
|\l(t)-\l(t+\r)|\leq \frac{\r\log n}{N(t+\r)}.
\eeq
Because $\l$ changes very little, simple estimates then give
\begin{claim}\label{lem++}
\begin{align}
&|A(t)-A(t+\r)|=O\bfrac{\r\log n}{N(t+\r)} &|B(t)-B(t+\r)|=
O\bfrac{\r\log n}{N(t+\r)}\label{pp1}\\
&|C(t)-C(t+\r)|=O\bfrac{\r\log n}{N(t+\r)}
&|D(t)-D(t+\r)|=O\bfrac{\r\log n}{N(t+\r)}\label{pp2}
\end{align}
If $||\bu(t)-\hu(t)||_1\leq n^{8/9}$ then
\begin{align}
&|A(t)-\hA(t)|=O\bfrac{||\bu(t)-\hu(t)||_1}{N(t)}
 &|B(t)-\hB(t)|=O\bfrac{||\bu(t)-\hu(t)||_1}{N(t)}\label{pp1a}\\
&|C(t)-\hC(t)|=
O\bfrac{||\bu(t)-\hu(t)||_1}{N(t)}&|D(t)-\hD(t)|=O\bfrac{||\bu(t)-\hu(t)||_1}{N(t)}\label{pp2a}
\end{align}
\end{claim}
\proofstart
The first expressions in \eqref{pp1} and \eqref{pp2} are easy to
deal with as the functions $f_j$ are smooth
and $\l$ is bounded throughout, see Lemma \ref{lambda}. Thus the each $f_j$ changes by
$O(\r\log n/N)$ and $y,z,\m$ change by $O(\r\log n)$ and $\m=\Omega(N)$.

For \eqref{pp1a} and \eqref{pp2a} we use Lemma \ref{onestep} to argue that
$$|\l(t)-\hl(t)|=O\bfrac{||\bu(t)-\hu(t)||_1}{N(t+\r)}.$$
Our assumption $t\leq T_0$ implies that $\m(t)=\Omega(n)$ and then
$\m(t)\sim \hm(t)$ and $N(t)\sim \hat{N}(t)$ and we can argue as for \eqref{pp1} and \eqref{pp2}.\\
{\bf End of proof of Claim \ref{lem++}}

Now fix $t$ and define for $\xi=y_1,y_2,z_1$,
$$
X_i(\xi)=\begin{cases}\xi(t+i+1)-\xi(t+i)&t+i<T_0\\ \ex[\xi(t+1)-\xi(t)\mid \bv(t)]&t+i\geq T_0
\end{cases}
$$
Then,
\beq{pp3}
\log n\geq\ex[X_i(\xi)\mid \bv(t+i)]=\sum_{\eta\in\set{a,b,c,2}}\d_\eta(t+i)
\psi_{\eta,\xi}(\bu(t+i))+O\bfrac{\log^2N(t+i)}{\l(t+i)N(t+i)}
\eeq
It follows from \eqref{pp0} -- \eqref{pp2} that
for all $\eta,\xi$ and $i\leq \r$,
$$\psi_{\eta,\xi}(\bu(t+i))=\psi_{\eta,\xi}(\bu(t))+O\bfrac{\log n}{n^{1-\a}}.$$
It then follows from \eqref{pp3} that \qs
\beq{pp4}
\log n\geq
\ex[\xi(t+\r)-\xi(t)\mid\bu(t)]=\sum_{\eta\in\set{a,b,c,2}}\D_\eta(\bu(t,t+\r))
\psi_{\eta,\xi}(\bu(t))+O\bfrac{\log n}{n^{1-\a}}
\eeq
This can be written as follows: We let
$\D_a= \D_a(\bu(t,t+\r))/\r$ etc. and $A=A(t),B=B(t)$.
\beq{equations1}
\begin{array}{rrrrl}
-\D_a&&&&=O(\r^{-1}\log n)\\
&(1-A)\D_b&&&=A+O(\r^{-1}\log n)\\
&-B\D_b&+\D_c&&=B+O(\r^{-1}\log n)\\
\D_a&+\D_b&+\D_c&+\D_2&=1.
\end{array}
\eeq
In comparison with \eqref{equations} we see, using \eqref{pp1a}, \eqref{pp2a} that
\beq{xi}
|\r\th_\xi(\hu(t))-\D_\xi|=O\brac{\log n+\frac{\r||\bu(t)-\hu(t)||_1}{N}}
\text{ for }\xi=a,b,c,2.
\eeq
Note that $A,\hA\leq 1/2$, see \eqref{2A}. This will be useful in dealing with
$\th_b$ and $\D_b$.

We now consider the difference between $\hu$ and $\bu$ at times $\r,2\r,\ldots$. We write
\beq{z1}
\xi(i\r)-\hx(i\r)=\xi((i-1)\r)-\hx((i-1)\r)+\sum_{t=(i-1)\r+1}^{i\r}([\xi(t)-\xi(t-1)]-[\hx(t)-\hx(t-1)])
\eeq
where $\xi=y,z,\m$ and $\hx=\hy,\hz,\hm$ in turn.
Then we write
\beq{z2}
\xi(t)-\xi(t-1)=\a_t+\b_t\text{ and }\hx(t)-\hx(t-1)=\ha_t+\hb_t
\eeq
where
$$\a_t=\sum_{\eta\in\set{a,b,c,2}}\d_{\eta,\xi}(\bu(t-1))\psi_{\eta,\xi}(\bu(t-1))
\text{ and }\b_t=\xi(t)-\xi(t-1)-\a_t$$
and
$$\ha_t=\sum_{\eta\in\set{a,b,c,2}}
\th_{\eta,\hx}(\hu(t-1))\psi_{\eta,\hx}(\hu(t-1))\text{ and }\hb_t=\hx(t)-\hx(t-1)-\ha_t.$$
It follows from \eqref{06q1}, \eqref{07q1} etc. that
$$\ex[\b_t\mid \bu(t-1)]=O\bfrac{\log^2N(t)}{\l(t)N(t)}.$$
An easy bound, which is a consequence of the Azuma-Hoeffding inequality, is that
\beq{z3}
\Pr\brac{\sum_{t=(i-1)\r+1}^{i\r}\b_t\ge \r^{1/2}\log^2n}\leq e^{-\Omega(\log^2n)}.
\eeq
We see furthermore that
\begin{align}
&\sum_{t=(i-1)\r+1}^{i\r}\hb_t=\nonumber\\
&=\sum_{t=(i-1)\r+1}^{i\r}\brac{\hx'(\hu(t-1+\varsigma_t))-\sum_{\eta\in\set{a,b,c,2}}
\th_{\eta,\hx}(\hu(t-1))\hx_{\eta}'(\hu(t-1))}\nonumber\\
&=\sum_{t=(i-1)\r+1}^{i\r}\brac{\hx'(\hu(t-1))+O\bfrac{\log n}{N}-\sum_{\eta\in\set{a,b,c,2}}
\th_{\eta,\hx}(\hu(t-1))\hx_{\eta}'(\hu(t-1))}\nonumber\\
&=O\bfrac{\r\log n}{N}=o(\r^{1/2}\log^2n),\label{z4}
\end{align}
where $0\leq\varsigma_t\leq1$ and $\hx'_\eta(t)$ is the derivative of $\hx$ in Case $\eta$.

In this and the following claims we take $N=N(i\r)$, the number of vertices at time $i\r$.

Now write
\begin{multline}\label{z5}
\sum_{t=(i-1)\r+1}^{i\r}\a_t=
\sum_{t=(i-1)\r+1}^{i\r}
\sum_{\eta\in\set{a,b,c,2}}\d_{\eta,\xi}(\bu(t))
\brac{\psi_{\eta,\xi}(\bu((i-1)\r)+O\bfrac{\r\log n}{N}}\\
=\sum_{\eta\in\set{a,b,c,2}}\D_\xi(\bu((i-1)\r+1,i\r))\psi_{\eta,\xi}(\bu((i-1)\r)
+O\bfrac{\r^2\log n}{N}
\end{multline}
and
\begin{multline}\label{z6}
\sum_{t=(i-1)\r+1}^{i\r}\ha_t=\\
\sum_{t=(i-1)\r+1}^{i\r}
\sum_{\eta\in\set{a,b,c,2}}\brac{\th_{\eta,\xi}(\hu((i-1)\r))+O\bfrac{\r}{N}}
\brac{\psi_{\eta,\xi}(\hu((i-1)\r)+O\bfrac{\r}{N}}\\
=\sum_{\eta\in\set{a,b,c,2}}\r\th_{\eta,\xi}(\hu((i-1)\r))\psi_{\eta,\xi}(\hu((i-1)\r)
+O\bfrac{\r^2}{N}.
\end{multline}
It follows that
\beq{z7}
\sum_{t=(i-1)\r+1}^{i\r}(\ha_t-\a_t)=A_1+A_2+o\brac{\r^{1/2}\log^2 n}
\eeq
where
\begin{align}
A_1&= \sum_{\eta\in\set{a,b,c,2}}(\D_\xi(\bu((i-1)\r+1,i\r))-\r\th_{\eta,\xi}(\hu((i-1)\r)))
\psi_{\eta,\xi}(\bu((i-1)\r)\nonumber\\
&=O\brac{\log n+\frac{\r||\bu((i-1)\r)-\hu((i-1)\r)||_1}{N}}.
\label{z8}\\
A_2&= \r\sum_{\eta\in\set{a,b,c,2}}\psi_{\eta,\xi}(\bu((i-1)\r)
(\psi_{\eta,\xi}(\bu((i-1)\r)-\psi_{\eta,\xi}(\hu((i-1)\r))\nonumber\\
&=O\brac{\frac{\r||\bu((i-1)\r)-\hu((i-1)\r)||_1}{N}}\label{z9}.
\end{align}
It follows from \eqref{z1} to \eqref{z9} that \whp, $i\r\leq T_0$ implies that with
\beq{ai=}
a_i=||\bu(i\r)-\hu(i\r)||_1
\eeq
that for some $C_1>0$,
$$a_i\leq a_{i-1}\brac{1+\frac{C_1\r}{N_i}}+2\r^{1/2}\log^2n
$$
where $N_i=N(i\r)\geq n/2$.

Putting
$$\Pi_i=\prod_{j=0}^i\brac{1+\frac{C_1\r}{N_j}}\leq e^{2C_1i\r/n}$$
we see by induction that
\beq{ai}
a_i\leq 2\r^{1/2}\log^2n\sum_{j=0}^i\frac{\Pi_i}{\Pi_j}\leq 2\r^{1/2}\log^2n(i+1)e^{2C_1i\r/n}.
\eeq
Since $i\leq n/\r$ we have
$$||\bu(i\r)-\hu(i\r)||_1=O(n\r^{-1/2}\log^2n).$$
Going from $\r\rdown{T_0/\r}$ to $T_0$ adds at most $\r\log n$ to the gap and the lemma follows.
\proofend

\section{Approximate equations}\label{approxeq}
The equations \eqref{slide} are rather complicated and we have not
made much progress in solving them. Nevertheless,
we can obtain information about them from a simpler set of equations
that closely approximate them when
$c$ is sufficiently large. The important observation is that when $\hl$ is large,
\beq{approx}
\hA\ll 1;\quad \hB\ll 1;\quad \hC\approx \frac{\hy\hl}{2\hm};
\quad \hD\approx \frac{\hz\hl^2}{2\hm};\quad \hl\approx \frac{2\hm}{\hy+\hz}.
\eeq

We will therefore approximate equations \eqref{slide} by the following
equations in variables $\tiy,\tiz,\tim$, $\til$:
\begin{align}
&\tiy'=-\frac{\tiy}{\tiy+\tiz}-1\label{hydash}\\
&\tiz'=\frac{2\tiy}{\tiy+\tiz}\label{hzdash}\\
&\tim'=-1-\frac{2\tiz\tim}{(\tiy+\tiz)^2}\label{hmdash}\\
&\til=\frac{2\tim}{\tiy+\tiz}.\label{hldash}
\end{align}
The initial conditions for $\tiy,\tiz,\tim,\til$ are that they
start out equal to $\hy,\hz,\hm,\hl$ at time $t=0$ i.e.
\beq{initial}
\tiy(0)=n;\quad \tiz(0)=0;\quad \tim(0)=cn;\quad \til=2c.
\eeq
\subsubsection{Analysis of the approximate equations}
The first two approximate equations imply $(\tiy +\tiz/2)^\prime=-1$, so that
$$
\tiy +\frac{\tiz}{2}=n-t.
$$
Using the second approximate equation and $\tiy=1-t-\tiz/2$, we obtain
$$
\tiz^\prime=\frac{2(n-t-\tiz/2)}{n-t+\tiz/2},
$$
or, introducing $\tau=n-t$ and
$$
X=\frac{\tiz}{2(n-t)}=\frac{\tiz}{2\tau},
$$
we get
\begin{equation}\label{dXdtau}
\frac{X+1}{X^2+1}\,dX =-\frac{1}{\tau}\,d\tau.
\end{equation}
Integrating,
$$
\frac{1}{2}\ln(X^2+1)+\arctan X =-\ln \tau +C.
$$
Now, at $t=0$ we have $\tau=n$ and $X=0$. So $C=\ln n$, i. e.
$$
\frac{1}{2}\ln(X^2+1)+\arctan X =-\ln(\tau/n).
$$
Let $\tiT$ satisfy $\tiy(\tiT)=0$.
At $t=\tiT$, we have $X=1$, so
$$
\ln n-\ln(n-\tiT)=\frac{1}{2}\ln 2 +\frac{\pi}{4}
$$
which implies
\beq{TT}
\tiT = \brac{1-\frac{1}{2^{1/2}}e^{-\pi/4}}n\approx 0.677603n.
\eeq
Note that
\begin{align*}
\til'&=\frac{2\tim'}{\tiy+\tiz}-\frac{2\tim(\tiy'+\tiz')}{(\tiy+\tiz)^2}\nonumber\\
&=-\frac{2}{\tiy+\tiz}-\frac{4\tiz\tim}{(\tiy+\tiz)^3}-\frac{2\tim}{(\tiy+\tiz)^2}
\brac{\frac{\tiy}{\tiy+\tiz}-1}
\nonumber\\
&=-\frac {2}{\tiy+\tiz}-\frac{2\tiz\tim}{(\tiy+\tiz)^3}\nonumber\\
&=-\frac {2}{\tiy+\tiz}-\frac{\tiz\til}{(\tiy+\tiz)^2},
\end{align*}
\beq{ldash}
\text{which implies that $\til$ is decreasing with $t$, at least as long as $\tiy,\tiz,\til>0$.}
\eeq
Here
\begin{align*}
\frac{\tiz}{(\tiy+\tiz)^2}=&\,\frac{\tiz}{(n-t+\tiz/2)^2}\\
=&\frac{\tiz}{(n-t)^2(1+X)^2}\\
=&\frac{2X}{(n-t)(1+X)^2}.
\end{align*}
Likewise
$$
-\frac{2}{\tiy+\tiz}=-\frac{2}{(n-t)(1+X)}.
$$
So $\til$ satisfies
$$
\til^\prime=-\frac{2}{(n-t)(1+X)}-\frac{2X}{(n-t)(1+X)^2}\,\til,\quad\til(0)=2c.
$$
Using \eqref{dXdtau}, we obtain
$$
\frac{d\til}{dX}=-\frac{2}{1+X^2}-\frac{2X}{(1+X)(1+X^2)}\,\til,\quad \left.\til(X)\right|_{X=0}=2c.
$$
Integrating this first-order, linear ODE, we obtain
\beq{comp}
\til(X)=\frac{(1+X)e^{-\arctan X}}{\sqrt{1+X^2}}
\left[2c-\int_0^X \frac{2e^{\arctan x}}{(1+x) \sqrt{1+x^2}}\,dx\right].
\eeq
In which case
\beq{iwc}
\til(\tiT)\approx 1.53c-1.418.
\eeq
\subsubsection{Simple Inequalities}\label{simpineq}
We will use the following to quantify \eqref{approx}:
$$
1\leq \frac{f_2(\hl)}{f_3(\hl)}= 1+\e_1,\qquad 1\le \frac{f_0(\hl)}{f_2(\hl)}=1+\e_2, \qquad
1\leq \frac{f_0(\hl)}{f_3(\hl)}=1+\e_3.$$
where
$$\e_1=\frac{\hl^2}{2f_3(\hl)},\qquad\e_2=\frac{1+\hl}{f_2(\hl)},\qquad\e_3=\frac{\hl^2+2\hl+2}{2f_3(\hl)}.$$
We use the above to verify the following sequence of inequalities
for $\hy,\hz,\hm,\hl$:
\begin{align}
&\frac{2\hm}{\hy+\hz}(1-\e_4)\leq\hl\leq \frac{2\hm}{\hy+\hz}.\label{first}\\
&0\leq \hA\leq \e_5\nonumber\\
&0\leq \hB\leq \e_6\nonumber\\
&\frac{\hy\hl}{2\hm}\leq \hC= \frac{y\hl}{2\hm}(1+\e_1).\nonumber\\
&\frac{z\hl^2}{2\hm}\leq \hD=\frac{\hz\hl^2}{2\hm}(1+\e_2).\nonumber
\end{align}
where
$$\e_4=\frac{\e_1}{1+\e_1},\qquad
\e_5=\frac{(1+\e_2)(1+\e_3)\l^3}{8f_0(\hl)},\qquad\e_6=\frac{\hl^2(1+\e_2)^2}{f_0(\hl)}.$$
(We use \eqref{slide2} to get $\hy\hz\leq \hm^2/\hl^2$ for use in defining $\e_5$).

For \eqref{first} we use
$$\frac{\hl(\hy+\hz)}{2\hm}\geq
\frac{1}{\max\set{\frac{f_2(\hl)}{f_3(\hl)},\frac{f_1(\hl)}{{f_2(\hl)}}}}=
\frac{f_3(\hl)}{f_2(\hl)}=
\frac{1}{1+\e_1}.$$

It follows from \eqref{first} that the initial value $\hl_0$ of $\hl$
satisfies
$$2c\geq \hl_0\geq 2c(1-\e_4).$$
Now $\e_4\leq .0001$ for $x\geq 15$ and so
\beq{lh0}
2c(1-.0001)\leq \hl_0\leq 2c.
\eeq
\subsubsection{Main Goal}\label{maingoal}
Lemma \ref{near} (below) in conjunction with Lemma \ref{close},
will enable us to argue that \whp\ in the process \2G, at some time
$T\leq T_0$ we will have
\beq{Texists}
y(T)=0,z(T)=\Omega(n)\text{ and }\l(t)=\Omega(1)\text{ for }t\leq T.
\eeq
Define
$$T_+=\min\set{t>0:\hy(t)\leq0\text{ or }\hz(t)\leq0\text{ or }
\tiy(t)\leq0\text{ or }\tiz(t)\leq0}.$$
We can bound this from below by a small constant as follows:
Initially $\hA,\hB$ are small $\hC$ is close to one for $c\geq 15$ and so \eqref{slide} implies
that $\hz$ is strictly increasing at the beginning. Also,
$\hy,\tiy$ start out large ($=n$) and so remain positive
initially.

Next define
\beq{T1def}
T_1=\min\set{T_+,\max\set{t:\hl(\t))\geq \l^*
\text{ and }\min\set{\tiy(\t)+\tiz(\t),\hy(\t)+\hz(\t)}\geq \b n\text{ for }
\t\leq t}}
\eeq
where
$$\b=-.01+\tiz(\tiT)/n=-.01+ 2(n-\tiT)/n\approx .63$$
and
$$\l^*=\til(\tiT)-5.$$
Comparing \eqref{iwc} and \eqref{lh0} we see that $\hl_0>\l^*$.

Note that
\beq{T0TT}
T_1\leq \tiT.
\eeq
This is because $\tiy(\tiT)=0$ and $\tiy'(\tiT)=-1$.
\begin{lemma}\label{near}
For large enough $c$,
\beq{hT}
\hy(T_1)=0<\hz(T_1)=\Omega(n)\text{ and }\hl(T_1)=\Omega(1).
\eeq
\end{lemma}
\proofstart
It follows from \eqref{slide} and Section \ref{simpineq}  that
\begin{align}
\hy'&\leq \frac{\e_6-\frac{\hy\hl}{2\hm}}{1-\e_5}-1\leq -\frac{\hy}{\hy+\hz}-1+\e_7\nonumber\\
\hy'&\geq -\frac{\frac{\hy\hl}{2\hm}(1+\e_1)}{1-\e_5}-1
\geq -\frac{\hy}{\hy+\hz}-1-\e_8.\label{lowerydash}\\ \nonumber\\
\hz'&\leq \frac{2\hy\hl}{2\hm}\cdot\frac{1+\e_1}{1-\e_5}
\leq \frac{2\hy}{\hy+\hz}+2\e_{8}.\nonumber\\
\hz'&\geq \frac{2\hy}{\hy+\hz}-2\e_{7}
\label{lowerzdash}\\\nonumber\\
\hl&\leq \frac{2\hm}{\hy+\hz}\label{lowerldash}\\
\hl&\geq \frac{2\hm}{\hy+\hz}(1-\e_4)
\geq \frac{2\hm}{\hy+\hz}-\e_{9}\label{upperldash}\\ \nonumber\\
\hm'&\leq -1-\frac{\hz\hl^2}{2\hm}
\leq -1-\frac{2\hz\hm}{(\hy+\hz)^2}(1-\e_4)^2\leq-1-\frac{2\hz\hm}{(\hy+\hz)^2}
+\e_{10}
\nonumber\\
\hm'&\geq -\frac{1+\frac{\hz\hl^2}{2\hm}(1+\e_2)}{1-\e_5}\geq
-1-\frac{2\hz\hm}{(\hy+\hz)^2}
-\e_{11}.\label{lowermdash}
\end{align}
where
\begin{align*}
&\e_7=\frac{\e_4+\e_5+\e_6}{1-\e_5},&\e_8=\frac{\e_1+\e_5}{1-\e_5},\qquad\qquad\qquad
&\e_{9}=\hl\e_4,\\
&\e_{10}=\frac{2\hl\e_4}{1-\e_4},&\e_{11}=\frac{\hl(\e_2+\e_5)+\e_5}{(1-\e_4)(1-\e_5)}
\end{align*}

When $t=0$ we have $\hy=n,\hz=0,\hm=cn$ and $\hl$
satisfying \eqref{lh0}, we see that $T_1>0$ for $c\geq 15$.

We can write $\hy(0)=n,\hz(0)=0,\hm(0)=cn$ and
\begin{align}
&\hy'=-\frac{\hy}{\hy+\hz}-1+\th_1\qquad&\text{ where }|\th_1|\leq \sd.\label{delta1}\\
&\hz'=\frac{2\hy}{\hy+\hz}+\th_2\qquad&\text{ where }|\th_2|\leq 2\sd.\label{delta2}\\
&\hm'=-1-\frac{2\hz\hm}{(\hy+\hz)^2}+\th_3&
\qquad\text{ where }-\e_{11}\le\th_3\le\e_{10}\label{delta3}\\
&\hl=\frac{2\hm}{\hy+\hz}+\th_4&\qquad\text{ where }
-\e_9\leq \th_4\leq 0.\label{delta4}
\end{align}
where
$$\sd=\max\set{\e_1,\e_2,\ldots,\e_8}.$$

It can easily be checked that the functions $\e_1,\ldots,\e_{11}$ are all
monotone decreasing for $\hl\geq \l^*$, ($\l^*(15)\approx 16.549$).
Furthermore, $\sd(16)<.00011$ and our error
estimates will mostly be $\sd$ times a moderate size constant. The only exceptions contain
a factor $c$, but if $c$ is large then $\sd$ will decrease to compensate.

It follows from \eqref{hmdash} and \eqref{delta3} that
\beq{mudown}
\hm,\tim\text{ both decrease for }t\leq T_1,\text{ since }\th_3<1\text{ for }\hl\geq \l^*.
\eeq

The ensuing calculations involve many constants and the expressions
\eqref{yzok} and \eqref{very} claim some inequalities that are tedious to justify.
It is unrealistic to expect the reader to check these calculations. Instead,
we have provided mathematica output in an appendix that will be seen to justify our claims.

The reader will notice the similarity between these equations and the approximation
\eqref{hydash} -- \eqref{hldash}. We will now refer to the equations \eqref{delta1} --
\eqref{delta4} as the {\em true} equations and \eqref{hydash} -- \eqref{hldash} as the
{\em approximate} equations.

\subsubsection{$y,z$ and $\hy,\hz$ are close}
We claim next that
\beq{half}
\max\set{|\hy(t)-\tiy(t)|,|\hz(t)-\tiz(t)|}\leq\sd F_1(t/n)n
\text{ for }0\leq t\leq T_1.
\eeq
where
\beq{fig}
F_a(x)=\b(e^{2ax/\b}-1)\text{ for }x\leq \frac{\tiT}{n}.
\eeq
for $a>0$.

 Note that
$$F_a'(t)=2(aF_a(t)/\b+1).$$
In the proof of \eqref{half}, think of $n$ as fixed and $h$ as a parameter that tends to zero.
Think of $\e$ as small, but fixed until the end of
the proof. In the display beginning with equation
\eqref{dis1}, only $h$ is the quantity going to zero.
Let
$$\hhu_i=\hy(ih),\hhv_i=\hz(ih),\tiu_i=\tiy(ih),\tiv_i=\tiz(ih)\text{ for }0\leq i\leq n/h.$$
Assume inductively that for $i<i_0=T_1/h$
\begin{equation}\label{ui}
|\hhu_i-\tiu_i|,|\hhv_i-\tiv_i|\leq \d F_{1+\e}(ih/n)n.
\end{equation}
This is true for $i=0$.

Suppose that
$$F_{1+\e}((i+1)h/n)=F_{1+\e}(ih/n)+\frac{h}{n}F_{1+\e}'((i+\th)h/n)$$
for some $0\leq \th\leq 1$.

Then by the inductive assumption and the Taylor expansion and uniform boundedness of second
derivatives,
\begin{align}
\hhv_{i+1}-\tiv_{i+1}&=\hhv_i-\tiv_i+h\brac{\frac{2\hhu_i}{\hhu_i+\hhv_i}-
\frac{2\tiu_i}{\tiu_i+\tiv_i}+\th_2(ih)}+O(h^2)\label{dis1}\\
&=\hhv_i-\tiv_i+h\brac{\frac{2\hhu_i(\tiv_i-\hhv_i)-
2\hhv_i(\tiu_i-\hhu_i)}{(\hhu_i+\hhv_i)(\tiu_i+\tiv_i)}+\th_2(ih)}+O(h^2)\nonumber\\
&\leq\hhv_i-\tiv_i+h\brac{\frac{2(\tiu_i+\tiv_i)
\max\set{|\hhv_i-\tiv_i|,|\hhu_i-\tiu_i|}}{(\hhu_i+\hhv_i)(\tiu_i+\tiv_i)}
+\th_2(ih)}+O(h^2)\nonumber\\
&\leq \sd F_{1+\e}(ih/n)n+2\sd h\brac{F_1(ih/n)/\b+1}+O(h^2)\nonumber\\
&=\sd(F_{1+\e}(ih/n)n+hF_1'(ih/n))+O(h^2)\nonumber\\
&=\sd(F_{1+\e}((i+1)h/n)n+h(F_1'(ih/n)-F_{1+\e}'((i+\th)h/n)))+O(h^2)\nonumber\\
&\leq \sd(F_{1+\e}((i+1)h/n)n-\Omega(\e h),\nonumber
\end{align}
completing the induction.

The remaining three cases are proved similarly.
This completes the inductive proof of \eqref{ui}.
Letting $\e\to 0$ we see for example that $\hy(t)-\tiy(t)\leq \sd F_1(t)$ for $t\leq T_1$.
This completes the proof of \eqref{half}.

Let
$$\a_0=F_1(\tiT/n).$$
Observe next that
\beq{hhalf}
(\tiy+\tiz)'=\frac{\tiy}{\tiy+\tiz}-1\leq 0.
\eeq
So for $t\leq T_1$ we have
\beq{644}
\tiy+\tiz\geq \tiy(\tiT)+\tiz(\tiT)=\tiz(\tiT)=2(n-\tiT)=(\b+.01)n.
\eeq
Furthermore, putting $X=1$ and going back to \eqref{comp},
\beq{l1}
\til(\tiT)=\frac{(1+\tiT)e^{-\arctan \tiT}}{\sqrt{1+\tiT^2}}
\brac{2c-\int_{x=0}^1\frac{2e^{\arctan x}}{(1+x) \sqrt{1+x^2}}\,dx}=
\a_1c-\a_2.
\eeq
\subsubsection{Lower bounding $\hl$}
We now show that $\til-\hl$ is small.
We now use \eqref{delta3} and \eqref{half} to write for $t\leq T_1$,
\begin{align*}
&|\tim'-\hm'|\\
&\leq |\th_3|+\card{\frac{2\hm \hz((\hy+\hz)^2+4\sd F_1(t/n)(\hy+\hz)n+4{\sd}^2 F_1(t/n)^2n^2)
-2\tim (\hy+\hz)^2(\hz-\sd F_1(t/n)n)}{(\hy+\hz)^2(\tiy+\tiz)^2}}\\
&=|\th_3|+\card{\frac{2\hz(\hy+\hz)^2(\hm-\tim)+2\sd F_1(t/n)n(\hy+\hz)
(4\hm \hz+\tim (\hy+\hz))+8\hm\hz{\sd}^2F_1(t/n)^2n^2}
{(\hy+\hz)^2(\tiy+\tiz)^2}}.
\end{align*}
Now, using \eqref{mudown},
$$\frac{4\hm\hz+\tim(\hy+\hz)}{(\hy+\hz)(\tiy+\tiz)^2}
\leq \frac{4\hm+\tim}{(\tiy+\tiz)^2}\leq \frac{5c}{\b^2n}$$
and
$$\frac{8\hm\hz}{(\hy+\hz)^2(\tiy+\tiz)^2}\leq \frac{8c}{\b^3n^2}.
$$

So,
$$|\tim'-\hm'|\leq \frac{2|\hm-\tim|}{\b n}+\a_3\sd$$
where
$$\a_3=\frac{10\a_0c}{\b^2}+\frac{8c\a_0^2\sd}{\b^3}+\frac{(2c+1)\sd}{\b(1-\sd)^2},$$
where the third term is an upper bound for $\e_{10},\e_{11}$ and its validity rests on 
\eqref{lowerldash} and \eqref{mudown}, with which we bound $\hl\leq 2c/\b$.

Integrating, we get that if $\z=|\hm-\tim|$ then
$$
\z'-\frac{2\z}{\b n}\leq \a_3\sd$$
and so
$$|\tim'-\hm'|\leq \a_3\sd e^{2t/\b n}\int_{\tau=0}^t e^{-2\tau/\b n}d\tau=\a_3\sd \frac{\b n}{2}(e^{2t/\b n}-1)\leq \a_4\sd n.$$
for $t\leq T_1$,
where
$$\a_4=\a_0\a_3/2.$$

It then follows that as long as $t\leq T_1$,
\begin{align*}
\til-\hl&=-\th_4+\frac{2\tim(\hy+\hz)-2\hm(\tiy+\tiz)}{(\hy+\hz)(\tiy+\tiz)}\\
&\leq \e_9+\frac{2(\hm+\a_4\sd n)(\hy+\hz)-
2\hm(\hy+\hz-2\a_0\sd n)}{(\hy+\hz)(\tiy+\tiz)}\\
&\leq \e_9+\frac{2\a_4\sd}{\b}+\frac{4c\a_0\sd}{\b^2}.
\end{align*}
It follows from \eqref{ldash} that for $t\leq T_1$ we have
\beq{lf}
\hl(t)\geq \til(T_1)-\a_5\sd
\eeq
where
$$\a_5=\frac{2c}{\b}+\frac{2\a_4}{\b}+\frac{4c\a_0}{\b^2}.$$

We now argue that $\hy(T_1)=0$ and $\hl(T_1)\geq \l^*$. This proves the Lemma \ref{near},
since $\hy(T_1)+\hz(T_1)\geq \b n$.
Suppose then to the contrary that $\hy(T_1)>0$.
 Recall that $T_1\leq \tiT$ (see \eqref{T0TT})
and suppose first that $T_1<\tiT$. Now let
$$T_2=\min\set{T_1+\e n,(T_1+\tiT)/2}$$
where $0<\e<10^{-10}$ is such that
\beq{defeps}
\max\set{\t\in [T_1,T_2]: \e\max\set{|\hl'(\t)|,|\hy'(\t)|,|\hz'(\t)|}\leq 10^{-10}}.
\eeq
The existence of such an $\e$ follows by elementary propositions in real analysis.

We will argue that $\t\in [T_1,T_2]$ implies
$$\hl(\t)\geq \l^*\text{ and }\min\set{\hy(\t)+\hz(\t),\tiy(\t)+\tiz(\t)}\geq \b n
\text{ and }\min\set{\hy(\t),\hz(\t),\tiy(\t),\tiz(\t)}>0,$$
which contradicts the definition of $T_1$.

Fix $\t\in [T_1,T_2]$. Now $\t<\tiT$ implies that $\tiy(\t)>0$. Together with
\eqref{hzdash} we see that $\tiz$ increases for $t\leq \tiT$ and hence $\tiz(\t)>0$.
We have $\tiy'(t)\geq -2$ (see \eqref{hydash})
and $\tiz'(t)\geq 0$ for $t\leq T_2$ (see \eqref{hzdash})
and so for some $\tau_1,\tau_2\in [T_1,T_2]$
\begin{align*}
\tiy(\t)&=\tiy(T_1)+(\t-T_1)\tiy'(\tau_1)\geq \tiy(T_1)-2\e n.\\
\tiz(\t)&=\tiz(T_1)+(\t-T_1)\tiz'(\tau_2)\geq \tiz(T_1).
\end{align*}
It follows (using \eqref{hhalf})  that
\beq{zz1}
\tiy(\t)+\tiz(\t)\geq \tiy(T_0)+\tiz(T_0)-2\e n\geq \tiy(\tiT)+\tiz(\tiT)-2\e n>\b n.
\eeq

We have, for some $\t_3,\t_4\in [T_1,T_2]$,
\begin{align}
\hy(\t)+\hz(\t)&=\hy(T_1)+\hz(T_1)+(\t-T_1)(\hy(\t_3)+\hz(\t_4))'\nonumber\\
&\geq \tiy(T_1)+\tiz(T_1)-2F_1(\tiT/n)\sd n-2\times 10^{-10}n\nonumber\\
&\geq \brac{\b+.01-2\brac{\a_0\sd+10^{-10})}}n\nonumber\\
&\geq \b n.\label{yzok}
\end{align}

We now argue that $\hz(\t)>0$.
Equation \eqref{hzdash} shows that $\tiz$ is strictly increasing initially.
Also, if $\hl\geq \l^*$ then
$\th_3\leq 1/8$. From \eqref{delta2}
we see that $\hz$ is strictly increasing at least
until a time $\t_0$ when $\hy(\t_0)\leq \b\sd$. On the other hand, we see from
\eqref{yzok} that if $\hy(\t)\leq \b\sd$ then $\hz(\t)>0$.  So,
\beq{fcv}
\min\set{\hy(\t),\hz(\t),\tiy(\t),\tiz(\t)}>0.
\eeq

Now we write
\begin{align}
\hl(\t)&=\hl(T_1)+(\t-T_1)\hl'(\t_3)\nonumber\\
&\geq \til(T_1)-(\til(T_1)-\hl(T_1))-10^{-10},\qquad\text{ using \eqref{defeps}},\nonumber\\
&\geq \a_1c-\a_2-\a_5\sd-10^{-10}.\nonumber\\
&>\l^*.\label{very}
\end{align}
We must now deal with the case where $T_1=\tiT$. Here we can just use \eqref{half} to argue
that $\hz(T_1)>\tiz(T_1)-\a_0\sd n>0$ and $\hy(T_1)+\hz(T_1)>\tiy(T_1)+\tiz(T_1)-\a_0\sd n>
(\b+.01-\a_0\sd)n>\b n$ and $\hl(\tiT)\geq \til(\tiT)-\a_5\sd>\l^*$.

This completes the proof of Lemma \ref{near}.
\proofend

It follows from Lemma \ref{close} that \whp\ $y(T_1)\leq n^{8/9}$,
$z(T_1)\geq \b n-n^{8/9}$ and $\l(T_1)\geq \l^*$. We claim that \qs,
$y$ becomes zero within the next $\n=n^{9/10}$ steps of
\2G. Suppose not.
It follows from Lemma \ref{onestep} that $\l$ changes by $o(1)$ and by \eqref{maxdegree} that
$z$ changes by $o(n)$ during these $\n$ steps. Thus $T_1+\n\leq T_0$.
It follows from \eqref{eqx5} that \qs\ at least $\n\log^{-2}n$ of these steps will be of type
Step 2. But each such step reduces $y$ by at least one, contradiction.

This verifies \eqref{Texists}.

\section{The number of components in the output of the algorithm}\label{noc}
We will tighten our bound on $\z$ from Lemma \ref{Klem}.
\begin{lemma}\label{Klemm}
If $c\geq 15$ then for
every positive constant $K$ there exists a constant $c_2=c_2(K)$ such that
$$\Pr\brac{\exists 1\leq t\leq T_1:\;\z(t)> c_2\log n}\leq n^{-K}.$$
\end{lemma}
\proofstart
We now need to use a sharper inequality than \eqref{x4} to replace $L_1$ by what is claimed in the
statement of the lemma. This sharper inequality uses higher moments of the $X_t$'s and we can
estimate them now that we have the estimate of the maximum of $\z(t)$ given in \eqref{x4}.
So, we now have to estimate terms of the form
$$\Psi_j(\xi\mid \boldsymbol{\eta})=\ex[\ |(\xi'-\xi)-\ex[(\xi'-\xi)\mid
\boldsymbol{\eta}]|^j\ \mid \boldsymbol{\eta}).$$
for $\xi=y_1,y_2,z_2$, $2\leq j\le \log n$ and $\boldsymbol{\eta}=\bv$ or $\bb,\bd$.

We use the inequality
$$(a+b+c+d)^j\leq 4^j(|a|^j+|b|^j+|c|^j+|d|^j)$$
for $j\geq 1$.

We will also need to estimate, for $2\leq j\leq \log n$,
\begin{multline*}
\sum_{k\geq 2}\frac{k(k-1)^j\l^k}{k!}=\l^2\sum_{k\geq 0}\frac{(k+1)^{j-1}\l^k}{(k-2)!}< 2^j\l^2\sum_{k\geq 0}\frac{k^j\l^k}{k!}
=2^j\l^2\sum_{k\geq 0}\sum_{\ell=0}^j\genfrac{\{}{\}}{0in}{}{j}{\ell}\frac{(k)_\ell\l^k}{k!}\\
=2^j\l^2\sum_{\ell=0}^j\genfrac{\{}{\}}{0in}{}{j}{\ell}\l^\ell\sum_{k\geq \ell}\frac{\l^{k-\ell}}{(k-\ell)!}
\leq 2^j\l^{j+2}e^\l \sum_{\ell=0}^j\genfrac{\{}{\}}{0in}{}{j}{\ell}\leq 2^jj!\l^{j+2}e^\l.
\end{multline*}
Here $\genfrac{\{}{\}}{0in}{}{j}{\ell}$ is a Stirling number of the second kind and it is easy to verify by induction on $j$ that the Bell number
$\sum_{\ell=0}^j\genfrac{\{}{\}}{0in}{}{j}{\ell}\leq j!$.

{\bf  Step 1.\/} $y_1+y_2+z_1>0$.
\si
{\bf Step 1(a)\/}. $y_1>0$.
\begin{eqnarray}
\Psi_j(y_1\mid\bb,\bd)&\leq&4^j\brac{{\frac{y_1}{2\mu}+
\sum_{k\ge 2}\frac{kz_k}{2\mu}\,(k-1)^j\frac{y_1}{2\mu}}
+\sum_{k\ge 2}\frac{kz_k}{2\mu}\,(k-1)^j\frac{2y_2}{2\mu}+\e_{\ref{04x1}}}.\label{04x1}\\
\Psi_j(y_i\mid\bv)&=&O\brac{2^{3j}\l^je^\l j!\brac{\frac{\z}{N}+\frac{\log^2N}{\l N}}}.\label{04xq1}\\
\nonumber\\ \nonumber\\
\Psi(y_2\mid \bb,\bd)&\leq&4^j\brac{\frac{2y_2}{2\mu}+
\sum_{k\ge 2}\frac{kz_k}{2\mu}\,(k-1)^j\frac{2y_2}{2\mu}
+\sum_{k\ge 2}\frac{kz_k}{2\mu}\,(k-1)^j\frac{3y_3}{2\mu}+\e_{\ref{041}}}.\label{041}\\
\Psi_j(y_2\mid |\bv|)&=&O\brac{2^{3j}\l^je^\l j!\brac{\l^3+\frac{\z}{N}+
\frac{\log^2N}{\l N}}}.\label{04q1}\\
\nonumber\\ \nonumber\\
\Psi(z_1\mid \bb,\bd)&\le&4^j\brac{\frac{z_1}{2\mu}+
\sum_{k\ge 2}\frac{kz_k}{2\mu}\,(k-1)^j\frac{z_1}{2\mu}
+\sum_{k\ge 2}\frac{kz_k}{2\mu}\,(k-1)^j\frac{2z_2}{2\mu}+\e_{\ref{051}}}.\label{051}\\
\Psi(z_1\mid |\bv|)&=&O\brac{2^{3j}\l^je^\l j!\brac{\l^2+\frac{\z}{N}+\frac{\log^2N}{\l N}}}.
\label{05q1}\\
\end{eqnarray}
{\bf Step 1(b)\/}. $y_1=0,y_2>0$.
\begin{eqnarray}
\Psi_j(y_1\mid \bb,\bd)&\leq&4^j\brac{
2\sum_{k\ge 2}\frac{kz_k}{2\mu}\,(k-1)^j\frac{2y_2}{2\mu}+\e_{\ref{4x1}}}.\label{4x1}\\
\Psi_j(y_1\mid |\bv|)&=&O\brac{2^{3j}\l^je^\l j!\brac{\frac{\z}{N}+
\frac{\log^2N}{\l N}}}.\label{4xq1}\\
\nonumber\\ \nonumber\\
\Psi_j(y_2\mid \bb,\bd)&\leq&4^j\brac{\frac{2y_2}{2\mu}+
2\sum_{k\ge 2}\frac{kz_k}{2\mu}\,(k-1)^j\frac{2y_2}{2\mu}
+2\sum_{k\ge 2}\frac{kz_k}{2\mu}\,(k-1)^j\frac{3y_3}{2\mu}+\e_{\ref{41}}}.\label{41}\\
\Psi_j(y_2\mid  |\bv|)&=&O\brac{2^{3j}\l^je^\l j!\brac{\l^3+\frac{\z}{N}+
\frac{\log^2N}{\l N}}}.\label{4q1}\\
\nonumber\\ \nonumber\\
\Psi_j(z_1\mid \bb,\bd)&\leq&4^j\brac{\frac{z_1}{\mu}+
2\sum_{k\ge 2}\frac{kz_k}{2\mu}\,(k-1)^j\frac{z_1}{2\mu}
+2\sum_{k\ge 2}\frac{kz_k}{2\mu}\,(k-1)^j\frac{2z_2}{2\mu}+\e_{\ref{51}}}.\label{51}\\
\Psi_j(z_1\mid
|\bv|)&=&O\brac{2^{3j}\l^je^\l j!\brac{\l^2+\frac{\z}{N}+\frac{\log^2N}{\l N}}}.\label{5q1mmm}\\
\nonumber\\ \nonumber\\
\end{eqnarray}
{\bf  Step 1(c).\/} $y_1=y_2=0,z_1>0$.
\begin{eqnarray}
\Psi_j(y_1\mid \bb,\bd)&=&\e_{\ref{9x1}}.\label{9x1}\\
\Psi_j(y_1\mid |\bv|)&=&\e_{\ref{9xq1}}.\label{9xq1}\\
\nonumber\\ \nonumber\\
\Psi_j(y_2\mid \bb,\bd)&\leq&2^j\brac{\sum_{k\ge 2}\frac{kz_k}{2\mu}\,(k-1)^j\frac{3y_3}{2\mu}
+\e_{\ref{91}}}.\label{91}\\
\Psi_j(y_2\mid |\bv|)&=&O\brac{2^{3j}\l^je^\l j!
\brac{\l^3+\frac{\z}{N}+\frac{\log^2N}{\l N}}}.\label{9q1}\\
\nonumber\\ \nonumber\\
\Psi_j(z_1\mid \bb,\bd)&\leq&4^j\brac{\frac{z_1}{2\mu}+
\sum_{k\ge 2}\frac{kz_k}{2\mu}\,(k-1)^j\frac{z_1}{2\mu}
+\sum_{k\ge 2}\frac{kz_k}{2\mu}\,(k-1)^j\frac{2z_2}{2\mu}+\e_{\ref{101}}}.\label{101}\\
\Psi_j(z_1\mid |\bv|)&=&O\brac{2^{3j}\l^je^\l j!
\brac{\l^2+\frac{\z}{N}+\frac{\log^2N}{\l N}}}.\label{10q1}
\end{eqnarray}
Now let $\cE_t=\set{\z(\t)\leq \log^2n:1\leq\t\leq t}$. Then let
$$Y_i=\begin{cases}(\z(i+1)-\z(i))1(\cE_i)
&0\leq i\leq T_1\\-c_1/2&T_1<i\leq n\end{cases}$$
Then, \qs
$$Y_{s+1}+\ldots+Y_t=\z(t)-\z(s)\text{ for }0\leq s<t\leq T_1.$$
For some absolute constant $c_2$, and with $\th=\frac{c_1}{100ce^{3ce+1}(3ce)^3}$ and $i\leq L_1$,
\begin{multline*}
\ex[e^{\th Y_{s+i}}\mid Y_{s+1},\ldots,Y_{s+i-1}]
=\sum_{k=0}^\infty \th^k\ex\left[\frac{Y_{s+i}^k}{k!}\bigg| Y_{s+1},\ldots,Y_{s+i-1}\right]\\
\leq 1-\th c_1/2
+c_2\sum_{k=2}^\infty \th^k2^{3k}\l(i)^{k+3}e^{\l(i)}\leq e^{-\th c_1/3},
\end{multline*}
where we have used \eqref{eqx2} and we have used Lemma \ref{lambda} to bound $\l(i)$.

It follows that for $t-s\leq L_1$ and real $u>0$
$$\Pr(Y_{s+1}+\cdots+Y_t\geq u)\leq e^{ -\th(u+ c_1(t-s)/3)}
$$
Suppose now that there exists $\t\leq T_0$ such that $\z(\t)\geq L_2$. Now q.s. there exists
$t_1\leq \t\leq t_1+L_1$ such that $\z(t_1)=0$. But then putting $u=-\log n$ and $L_2=\frac{6K\log n}{c_1}$ we see that
given $t_1$,
$$\Pr(\exists t_1\leq \t\leq t_1+L_1:\z(\t)\geq L_2)\leq 
\Pr\brac{\neg\bigcup_t\cE_t}+e^{-\th(c_1L_2/3-\log n)}\leq n^{-K}.$$
\proofend

We get a new path for every increase in $V_{0,j},\,j\leq 1$. If we look at equations
\eqref{04x} etc., then we see that
the expected number added to $V_{0,j}$ at step $t$ is $O(\z(t)/\m(t))$.
So if
$Z_P(t)$ is the number of increases at time $t$ and $Z_P=\sum_{t=0}^{T_3}Z_P(t)$, where $T_3$ is the
time at the beginning of Step 3, then
\beq{ZP}
\ex[Z_P]=O\ex\brac{\brac{\log n\sum_{t=0}^{T_3}\frac{1}{\mu(t)}}
=O\brac{\log n\ \ex\left[\log\bfrac{\m(0)}{\m(T_3)}\right]}}.
\eeq
Now in our case $\m(T_3)=\Omega(n)$ with probability $1-o(n^{-2})$
in which case $\ex[Z_p]=O(\log n)$.
We will apply the Chebyshev inequality to show concentration around the mean. We will
condition on $||\bu(t)-\hu(t)||_1\leq n^{8/9}$ for $t\leq T_1$ (see Lemma \ref{close}).
With this conditioning, the expected value of $Z_P(t)$ is determined up to a factor
$1-O(n^{-1/9}\log^2n)$ by the value of $\hu(t)$. In which case, $\ex[Z_P(t)\mid Z_P(s)]
=(1+o(1))\ex[Z_P(t)]$ and we can apply the Chebychev inequality to show that \whp\
$Z_P=O(\log n)$. We combine this with Lemma \ref{few} to obtain Theorem \ref{th1}.

\section{Hamilton cycles}\label{posa}
We will now show how we can use Theorem \ref{th1}(a) to prove the existence and construction of
Hamilton cycles. We will first need to remove a few random edges $X$ from $G=G_{n,cn}^{\d\geq 3}$
in such a way that the pair $(G-X,X)$ is distributed very close
to $(H=G_{n,cn-|X|}^{\d\geq 3},Y)$ where
$Y$ is a random set of edges disjoint from $E(H)$. In which case
we can apply Theorem \ref{th1} to $H$ and then we can use the edges of $Y$ to close cycles in the
extension-rotation procedure.
\subsection{Removing a random set of edges}\label{remove}
Let
$$s=n^{1/2}\log^{-2}n$$
and let
$$\Omega=\set{(H,Y):H\in \cG_{n,cn-s}^{\d\geq 3},
Y\subseteq\binom{[n]}{2},|Y|=s\text{ and }E(H)\cap Y=\emptyset}$$
where $\cG_{n,m}^{\d\geq 3}=\set{G_{n,m}^{\d\geq 3}}$.

We consider two ways of randomly choosing an element of $\Omega$.
\begin{enumerate}[{\bf (a)}]
\item First choose $G$ uniformly from $\cG_{n,cn}^{\d\geq 3}$ and then choose
an $s$-set $X$ uniformly from $E(G)\setminus E_3(G)$,
where $E_3(G)$ is the set of edges of $G$ that are
incident with a vertex of degree 3. This produces a pair
$(G-X,X)$. We let $\Pr_a$ denote the induced probability
measure on $\Omega$.
\item Choose $H$ uniformly from $\cG_{n,cn-s}^{\d\geq 3}$ and
then choose an $s$-set $Y$ uniformly from $\binom{[n]}{2}\setminus E(H)$.
This produces a pair
$(H,Y)$. We let $\Pr_b$ denote the induced probability
measure on $\Omega$.
\end{enumerate}
The following lemma implies that as far as properties that
happen \whp\ in $G$, we can use Method (b), just as well as Method (a) to
generate our pair $(H,Y)$.
\begin{lemma}\label{contig}
There exists $\Omega_1\subseteq \Omega$ such that
\begin{enumerate}[{\bf (i)}]
\item $\Pr_a(\Omega_1)=1-o(1)$.
\item $\om=(H,Y)\in \Omega_1$ implies that $\Pr_a(\om)=(1+o(1))\Pr_b(\om)$.
\end{enumerate}
\end{lemma}
\proofstart
We first compute the expectation of the number $\m_3=\m_3(G)$ of edges
incident to a vertex of degree 3 in $G$ chosen uniformly
from $\cG_{n,cn}^{\d\geq 3}$. We will use the random sequence model of
Section \ref{mod}. We will show that
$\m_3$ is highly concentrated in this model and then we can transfer this
result to our graph model. Observe
first that if $\n_3$ is the number of vertices of degree 3 in $G_\bx$ then
Lemma \ref{lem4x} implies that
$$\card{\n_3-\frac{\l^3}{3!f_3(\l)}n}=O(n^{1/2}\log n),\qquad\qs.$$
Here $\l$ is the solution to $\l f_2(\l)/f_3(\l)=2cn$.

To see how many edges are incident to these $\n_3$ vertices we consider the following experiment:
Condition on $\n_3=\r n$ where $\r$ will be taken to be close to $\r_3=\frac{\l^3}{3!f_3(\l)}$.
We take a
random permutation $\p$ of $[2cn]$ and compute the number $Z$ of $i\leq cn$
such that $\set{\p(2i-1),\p(2i)}\cap [3\n_3]\neq
\emptyset$. This will give us the number of edges in $G_\bx$ that are incident
with a vertex of degree 3. Now
$$\ex[Z]=cn\brac{1-\frac{2cn-\r n}{2cn}\frac{2cn-\r n-2}{2cn-2}}=cn(2\r-\r^2+O(1/n)).$$
Now interchanging two positions in $\p$ can change $Z$ by at most one and so
applying the Azuma-Hoeffding inequality
for permutations (see for example Lemma 11 of Frieze and Pittel \cite{FP1} or
Section 3.2 of McDiarmid \cite{McD}) we see that
$\Pr(|Z-\ex[Z]|\geq u)\leq e^{-u^2/(cn)}$ for any $u\geq 0$. Putting this all together we see that
$$\Pr(|\m_3(G)-\r_3(2-\r_3)cn|\geq u)\leq e^{-u^2/cn}.$$
Now let
$$\widehat{\cG}_{n,cn}^{\d\geq 3}=\set{G\in\cG_{n,cn}^{\d\geq 3}:
|\m_3(G)-\r_3(2-\r_3)cn|\leq n^{1/2}\log n}$$
and
$$\Omega_a=\set{(H,Y)\in\Omega:H+Y\in \widehat{\cG}_{n,cn}^{\d\geq 3}}.$$
This satisfies requirement (a) of the lemma.

Suppose next that $\om\in \Omega_a$. Then
\begin{align}
&\Pr_a(\om)=\frac{1}{|\cG_{n,cn}^{\d\geq 3}|}
\cdot\frac{1}{\binom{cn(1-\r_3)^2\pm n^{1/2}\log n}{s}}=
\frac{1+O(\log^{-1} n)}{|\cG_{n,cn}^{\d\geq 3}|\cdot\binom{cn(1-\r_3)^2}{s}}\label{obo1}\\
&\Pr_b(\om)=\frac{1}{|\cG_{n,cn-s}^{\d\geq 3}|}\cdot\frac{1}{\binom{\binom{n}{2}-cn}{s}}\label{obo2}
\end{align}
One can see from this that one has to estimate the ratio
$|\cG_{n,cn}^{\d\geq 3}|/ |\cG_{n,cn-s}^{\d\geq 3}|$.
For this we make estimates of
$$M=|\set{(G_1,G_2)\in\cG_{n,cn}^{\d\geq 3}\times
\cG_{n,cn-s}^{\d\geq 3}:E(G_1)\supseteq E(G_2) }|.$$
We have the following inequalities:
\begin{align}
&|\widehat{\cG}_{n,cn}^{\d\geq 3}|\binom{cn(1-\r_3)^2-n^{1/2}\log n}{s}\leq
M\leq|\widehat{\cG}_{n,cn}^{\d\geq 3}| \binom{cn(1-\r_3)^2+n^{1/2}\log n}{s}+\nonumber\\
&\gap{3}|\cG_{n,cn}^{\d\geq 3}|\sum_{|u|\geq n^{1/2}\log n}\binom{cn(1-\r_3)^2+u}{s}e^{-u^2/cn}\label{obo3}\\
&M=|\cG_{n,cn-s}^{\d\geq 3}|\binom{\binom{n}{2}-cn}{s}. \label{obo4}
\end{align}
We get \eqref{obo3} by summing $\m_3(G_1)$ over $G_1\in\cG_{n,cn}^{\d\geq 3}$ and bounding
$\m_3(G_1)$ according to whether or not $G$ is in $\widehat{\cG}_{n,cn}^{\d\geq 3}$.
Equation \eqref{obo4} is obtained by summing over $G_2\in \cG_{n,cn-s}^{\d\geq 3}$, the
number of ways of adding $s$ edges to $G_2$.

Now
\begin{multline*}
\sum_{|u|\geq n^{1/2}\log n}\binom{cn(1-\r_3)^2+u}{s}e^{-u^2/cn}\leq 2\sum_{u\geq n^{1/2}\log n}\binom{cn(1-\r_3)^2}{s}e^{O(us/n)}e^{-u^2/cn}\\
2\binom{cn(1-\r_3)^2}{s}\sum_{u\geq n^{1/2}\log n}e^{-u^2/2cn}=O\brac{\binom{cn(1-\r_3)^2}{s}e^{-\Omega(\log^2n)}}.
\end{multline*}

It follows from this and \eqref{obo3} that
$$M=|\cG_{n,cn}^{\d\geq 3}|\binom{cn(1-\r_3)^2}{s}\brac{1+O(\log^{-1}n)}.$$
By comparing with \eqref{obo4} we see that
$$\frac{|\cG_{n,cn}^{\d\geq 3}|}{|\cG_{n,cn-s}^{\d\geq 3}|}=(1+o(1))
\frac{\binom{\binom{n}{2}-cn}{s}}{\binom{cn(1-\r_3)^2}{s}}.$$
The lemma follows by using this in conjunction with \eqref{obo1} and \eqref{obo2}.
\proofend
\subsection{Connectivity of $G_{n,cn}^{\d\geq 3}$}\label{conn}
\begin{lemma}\label{conn1}
$G_{n,cn}^{\d\geq 3}$ is connected, \whp.
\end{lemma}
\proofstart
It follows from Lemma \ref{contig} that we can replace $G_{n,cn}^{\d\geq 3}$ by
$G_{n,cn-s}^{\d\geq 3}$ plus $s$ random edges.
We use the random sequence model to deal with
$G_{n,cn-s}^{\d\geq 3}$. Let
Fix $4\leq k\leq n/\log^{20}n$. For $K\subseteq [n]$,
$e(K)$ denotes the number of edges
of $G_\bx$ contained in $K$. Let $\ell_0=\log n/(\log\log n)^{1/2}$. Then
with $\l$ the solution to $\l f_2(\l)/f_3(\l)=2c$, 
\begin{align}
&\Pr(\exists K\subseteq [n]:e(K)\geq 5k/4)\leq o(1)+\d_k\binom{n}{k}\sum_{d=3k/2}^{\ell_0k}
\frac{\l^dk^d}{d!f_3(\l)^k}\binom{cn}{5k/4}\bfrac{d}{cn}^{5k/2}.\label{mm1}\\
&\leq \d_k\sum_{d=3k/2}^{\ell_0k}\bfrac{\l ek}{d}^d\bfrac{e^{9/4}\ell_0^{5/2}k^{1/4}}{(5c/4)^{5/4}f_3(\l)n^{1/4}}^k\leq 
\d_k\ell_0k\bfrac{e^{9/4}\ell_0^{5/2}k^{1/4}e^\l}{(5c/4)^{5/4}f_3(\l)n^{1/4}}^k\label{mm2}
\end{align}
{\bf Explanation of \eqref{mm1}:}
Here $\d_k=1+o(1)$ for $k\leq \log^2n$ and $O(n^{1/2})$ for larger $k$. The term
$\frac{\l^dk^d}{d!f_3(\l)^k}$ bounds the probability that the total degree of $K$ is $d$,
see \eqref{boundd}. Given the degree sequence we take a random permutation $\p$ of the multi-set
$\set{d_\bx(j)\times j:j\in[n]}$ and bound the probability that there is a set of $5k/4$ indices
$i$ such that $\p(2i-1),\p(2i)\in K$. This expression assumes that vertex degrees are
independent random variables. We can always inflate the estimate by $O(n^{1/2})$ to account
for the degree sum being fixed. This is what $\d_k$ does for $k\geq \log^2n$. For smaller $k$
we use \eqref{ll1}. The bound of $d\leq \ell_0k$ arises from Lemma \ref{lem4}(b).

Let $\s_{k}$ denote the RHS of \eqref{mm2}. Then, we have
$\sum_{k=4}^{n/\log^{20}n}\s_k=o(1)$.

But if no $G$ has minimum degree at least 3 and $K$ contains at most $5|K|/4$ edges then
there must be edges with one end in $K$. So, we see that \whp\ the minimum component size
in $G$ will be at least $n/\log^{20}n$. We now use the result of Section \ref{remove}.
If we take $H=G_{n,cn-s}^{\d\geq 3},\,s=n^{1/2}\log^{-2}n$ then we know by the above that \whp\ it
only has components of size at least $n/\log^{20}n$. Now add $s$ random edges $Y$. Then
$$\Pr(H+Y\text{ is not connected})=o(1)+\log^{40}n\brac{1-\frac{1}{\log^{40}n}}^s=o(1).$$
Now apply Lemma \ref{contig}.
\proofend
\subsection{Extension-Rotation Argument}\label{exro}
We will as in Section \ref{conn} replace $G_{n,cn}^{\d\geq 3}$ by
$G_{n,cn-s}^{\d\geq 3}$ plus $s$ random edges $Y$.
Having run \2G\ we will \whp\ have a two matching $M_0$ say such that $M_0$
has $O(\log n)$ components.

The main idea now of course is that of a {\em rotation}. Given a path $P=(u_1,u_2,
\ldots,u_k$ and an edge $(u_k,u_i)$ where $i\leq k-2$ we say that the path
$P'=(u_1,\ldots,u_i,u_k,u_{k-1},\ldots,u_{i+1})$ is obtained from $P$ by a rotation.
$u_1$ is the {\em fixed} endpoint of this rotation. We now describe an algorithm,
\HAM\ that \whp\ converts $M_0$ into a Hamilton cycle in $O(n^{1.5+o(1)})$ time.

Given a path $P$ with endpoints $a,b$ we define a {\em restricted rotation search}
$RRS(\n)$ as follows: Suppose that we have a path $P$ with endpoints $a,b$. We start
by doing a sequence of rotations with $a$ as the fixed endpoint. Furthermore
\begin{enumerate}[R1]
\item We only do a rotation if the endpoint of the path created is not an endpoint of the
paths that have been created so far.
\item We stop this process when we have either (i) created $\n$ endpoints or (ii)
we have found a path $Q$ with an endpoint that is outside $Q$. We say that we have found an
{\em extension}.
\end{enumerate}
Let $END(a)$ be the set of endpoints, other than $a$, produced by this procedure. The main
result of \cite{FP} is that \whp\ that regardless of our choice of path $P$,
either (i) we find an extension or (ii) we are able to
generate $n^{1-o(1)}$ endpoints. We will run this procedure with $\n=n^{3/4}\log^2n$.

Assuming that we did not find an extension and having constructed $END(a)$, we take each
$x\in END(a)$ in turn and starting with the path $P_x$
that we have found from $a$ to $x$, we carry out R1,R2 above with $x$ as the fixed
endpoint and either find an extension or create a set of $\n$ paths with
$x$ as one endpoint and the other endpoints comprising a set $END(x)$ of size $\n$.

It follows from \cite{AV} that the above construction $RSS(\n)$ can be carried out in
$O(\n^2\log n)$ time.

Algorithm \HAM
\begin{enumerate}[Step 1]
\item Choose a path component $P$ of the current 2-matching $M$, with endpoints $a,b$.\\
If there are no such components and $M$ is not a Hamilton cycle,
choose a cycle $C$ of $M$ and delete an edge to create $P$:
\item Carry out $RSS(\n)$ until either an extension is found or we have constructed
$\n+1$ endpoint sets.
\begin{description}
\item[Case a:]
We find an extension. Suppose that we construct a path $Q$ with endpoints $x,y$ such
that $y$ has a neighbour $z\notin Q$.
\begin{enumerate}[(i)]
\item If $z$ lies in a cycle $C$ then let $R$ be a path
obtained from $C$ by deleting one of the edges of $C$ incident with $z$. Let now
$P=x,Q,y,z,R$ and go to Step 1.
\item If $z=u_j$ lies on a path $R=(u_1,u_2,\ldots,u_k)$ where the numbering is chosen so
    that $j\geq k/2$ then we let $P=x,Q,y,z,u_{j-1},\ldots,u_1$ and go to Step 1.
\end{enumerate}
\item[Case b:] If there is no extension then we search for an edge $e=(p,q)\in Y$ such that
$p\in END(a)$ and $q\in END(p)$. if there is no such edge then the algorithm fails.
If there is such an edge, consider the cycle $P+e$.
Now either $C$ is a Hamilton cycle and we are done, or else there is a vertex $u\in C$ and
a vertex $v\notin C$ such that $(u,v)$ is an edge of $H$. Assuming that $H$ is connected,
see Lemma \ref{conn1}. We now delete one of the edges, $(u,w)$ say,
of $C$ incident with $u$ to create
a path $Q$ from $w$ to $u$ and treat $e$ as an extension of this path. We can now proceed as
Case a.
\end{description}
\end{enumerate}
\subsubsection{Analysis of \HAM}
We first bound the number of executions of $RSS(\n)$. Suppose that $M_0$ has $\k\leq K_1\log n$
components for some $K_1>0$.
Each time we execute Step 2, we either reduce the number of components by one or
we halve the size of one of the components not on the current path. So if the component
sizes of $M_0$ are $n_1,n_2,\ldots,n_\k$ then the number of executions of Step 2 can be
bounded by
$$\k+\sum_{i=1}^\k \log_2n_i\leq \k+\k\log_2(n/\k)=O(\log n\log\log n).$$
(Re-call that $\log(ab)\leq 2\log((a+b)/2)$ for $a,b>0$
and you will see the first inequality here).

An execution of Step 2 takes $O(\n^2\log n)$ time and so we are within the time bound claimed
by Theorem \ref{th2}.

We first argue that \HAM\ succeeds \whp. Suppose that the edges of $Y$ are $e_1,e_2,\ldots,e_s$.
We can allow the algorithm to access these edges in order, never going back to a previously
examined edge. The probably that an $e_i$ can be used in Case b is always at least
$\frac{\binom{\n}{2}-s}{\binom{n}{2}}\geq \frac{\log^4n}{2n^{1/2}}$ (we have subtracted $s$
because some of the useful edges might have been seen before the current edge in the order).
So the probability of failure is bounded by the probability that the binomial
$Bin\brac{s,\frac{\log^4n}{2n^{1/2}}}$ is less than $K_2\log n\log\log n$
for some $K_2>0$. And this tends to zero. This completes the proof of Theorem \ref{th2}.

\section{Concluding remarks}\label{final}
The main open question concerns what happens when $c<15$.
Is it true that \eqref{hT} holds all the way down
to $c>3/2$? We have done some numerical experiments and here are some results from these experiments:
$$
\begin{array}{ccccc}
c&y_{final}&z_{final}&\m_{final}&\l_{final}\\
3.0&0.000008&0.283721&0.398527&1.822428\\
2.9&0.000009&0.242563&0.326139&1.602749\\
2.8&0.000010&0.197461&0.253645&1.370798\\
2.7&0.000010&0.148901&0.182327&1.123928\\
2.6&0.000010&0.098344&0.114494&0.858355\\
2.5&0.000010&0.048976&0.054010&0.565840
\end{array}
$$
These are the results of running Euler's method with step length $10^{-5}$ on the sliding trajectory
\eqref{slide}. They indicate that \eqref{hT} holds
down to somewhere close to 2.5. This would indicate
some sort of phase transition in the performance of \2G\ at around this point. There is one for the
Karp-Sipser matching algorithm and so we are led to conjecture there is one here too.

Can we prove anything for $c<15$? At the moment we can not even show that at the completion of \2G\
the 2-matching $M$ has $o(n)$ components. This will be the
subject of further research.

Finally, we mention once again, the possible use of the ideas of \cite{CFM} to reduce the running
time of our Hamilton cycle algorithm to $O(n^{1+o(1)})$ time.

Our list of problems/conjectures arising from this research can thus be summarised:
\begin{enumerate}[{\bf (a)}]
\item Find a threshold $c_1$ such that \2G\ produces a 2-matching in $G_{n,cn}^{\d\geq 3}$ with
$O(\log n)$ components \whp\ iff $c>c_1$ .
\item If $c_1>3/2$ then show that when $c\in (3/2,c_1)$, the number of components
in the 2-matching produced is $O(n^\a)$ for some constant $\a<1$.
\item Analyse the performance of \2G\ on the random graph $G_{n,cn}$ i.e. do not condition
on degree at least three. Is there a threshold $c_2$ such that if $c\leq c_2$ then \whp\
only Steps 1a,1b,1c are needed, making the matching produced optimal.
\item Can \2G\ be used to find a Hamilton cycle \whp\ in $O(n^{1+o(1)})$ time when
applied to $G_{n,cn}^{\d\geq 3}$ and $c$ sufficiently large?
\item How much of this can be extended to find edge disjoint Hamilton cycles in
 $G_{n,cn}^{\d\geq k}$ for $k\geq 4$.
\end{enumerate}

\appendix
\section{Proof of \eqref{ll1}}
To find a sharp estimate for the probabilities in (\ref{ll1}) we have to refine a bit
the proof of the local limit theorem, since in our case the
variance of the $Z_ j$ are not always bounded away from zero. However it is enough to
consider the case where $N\s^2\to\infty$. There is little loss of generality
in assuming that $D=0$ here. As usual, we start
with the inversion formula
\begin{eqnarray}
\Pr\left(\sum_{j=1}^NZ_j=\t\right)&=&\frac
{1}{2\pi}\int_{-\pi}^{\pi}e^{-i\t x}
\E\left(e^{ix\sum_{j=1}^NZ_j}\right)\,dx\nonumber\\
&=&\frac {1}{2\pi}
\int_{-\pi}^{\pi}e^{-i\t x}\prod_{\ell=2}^3\left[\E(e^{ix\cP_\ell})\right]^{N_\ell}
\,dx,
\label{53}
\end{eqnarray}
where $\t=2M-k$.
Consider first $|x|\ge (N{\l})^{-5/12}$. Using an inequality (see Pittel
\cite{Pi})
$$
|f_\ell(\eta)|\le e^{({\rm Re}\eta -|\eta|)/(\ell+1)}f_\ell(|\eta|),
$$
we estimate
\begin{align}
&\frac {1}{2\pi}\int_{|x|\ge (N{\l})^{-5/12}}
\left|e^{-i\t x}\prod_{\ell=2}^3\left(\frac
{f_\ell(e^{ix}{\l})}{f_\ell({\l})}\right)^{N_\ell}\right|\,dx\nonumber\\
&\leq\frac {1}{2\pi}\int_{|x|\ge (N{\l})^{-
5/12}} e^{N{\l}(\cos x-1)/4}\,dx\nonumber\\
&\leq e^{N{\l}[(\cos ((N{\l})^{-5/12})-1)/4]}\nonumber\\
&\leq e^{-(N{\l})^{1/6}/9}. \label{54}
\end{align}
For $|x|\le (N{\l})^{-5/12}$, putting $\eta={\l}e^{ix}$ and using
$$\sum_{\ell=2}^3\frac{N_\ell {\l}f_\ell^\prime({\l})}{f_\ell({\l})}=2M\text{ and }
d/dx=i\eta d/d\eta$$
we expand $\sum_{\ell=2}^3N_\ell \log\brac{\frac{f_\ell(\eta)}{f_\ell({\l})}}$
as a Taylor series around $x=0$ to obtain
\begin{eqnarray}
-i\t x+\sum_{\ell=2}^3N_\ell\log \left(\frac {f_\ell(e^{ix}{\l})}{f_\ell({\l})}\right)&=&ikx-\frac {x^2}2
\left.{\cal D}
\brac{\sum_{\ell=2}^3N_\ell\frac {\eta f_\ell^\prime(\eta)}{f_\ell(\eta)}}\right|_{\eta={\l}}
\nonumber \\
&&-\frac {ix^3}{3!}\left.{\cal D}^2
\brac{\sum_{\ell=2}^3N_\ell\frac {\eta f_\ell^\prime(\eta)}{f_\ell(\eta)}}
\right|_{\eta={\l}}\nonumber\\
&&+O\brac{x^4\left.
{\cal D}^3\brac{\sum_{\ell=2}^3N_\ell\frac {\eta f_\ell^\prime(\eta)}{f_\ell(\eta)}}
\right|_{\eta=\tilde\eta}}.\qquad\,
\label{55}
\end{eqnarray}
Here $\tilde\eta={\l}e^{i\tilde x}$, with $\tilde x$ being
between $0$ and $x$, and ${\cal D}=\eta(d/d\eta)$. Now, the coefficients of
$x^2/2,\,
x^3/3!$ and $x^4$ are $N\s^2,\,O(N\s^2),\,O(N\s^2)$ respectively, and
$\s^2$ is of order ${\l}$. (Use (\ref{30}) and consider the effect of
${\cal D}$ on
a power of $\eta$.) So the second
and the third terms in (\ref{55}) are $o(1)$
uniformly for $|x|\le (N{\l})^{-5/12}$. Therefore
\begin{equation}
\frac 1{2\pi}\int_{|x|\le (N{\l})^{-5/12}}=\int_1+\int_2+\int_3, \label{55a}
\end{equation}
where
\begin{eqnarray}
\int_1&=&\frac {1}{2\pi}\int_{|x|\le (N{\l})^{-5/12}} e^{ikx-N\s^2x^2
/2}\,dx\nonumber\\
&=&\frac {1}{\sqrt {2\pi N\s^2}}+O\left(\frac{k^2+1}{({\l}N)^{3/2}}\right),\label{56}\\
\int_2&=&O\left({\cal D}^2
\brac{\sum_{\ell=2}^3N_\ell\frac {{\l} f_\ell^\prime({\l})}{f_\ell({\l})}}
\int_{|x|\le (N{\l})^{-5/12}} x^3e^{-N\s^2x^2/2}\,dx\right)\nonumber\\
&=&O\left(N{\l}\int_{ |x|\le (N{\l})^{-5/12}} |x|^3e^{-N\s^2x^2/2}\,
dx\right)\nonumber\\
&=&O(e^{-\alpha (N{\l})^{1/6}}),\label{57}\\
\noalign{($\alpha >0$ is an absolute constant), and}
\int_3&=&O\brac{N{\l}\int_{|x|\le (N{\l})^{-5/12}} x^4 e^{-N\s^2x^2/2}\,dx}
\nonumber\\
&=&o\left(\int_2\right).\label{58}
\end{eqnarray}
Using (\ref{53})-(\ref{58}), we arrive at
$$
\Pr\left(\sum_{\ell}Z_{\ell}=\t\right)=\frac {1}{\sqrt {2\pi N\s^2}}
\times\brac{1+O\left(\frac {k^2+1}{N{\l}}\right)}.
$$

\section{Mathematica Output}
In the computations below, $\e_p(\hl)$ is represented by $ep[x]$ and
$\a_p$ is represented by $Ap$ and $\b$ is represented by $B$. The computation $C_1$
is the justification for \eqref{very}.

\noindent\(\pmb{\text{f0}[\text{x$\_$}]\text{:=}\text{Exp}[x]}\)

\noindent\(\pmb{\text{f1}[\text{x$\_$}]\text{:=}\text{f0}[x]-1}\)

\noindent\(\pmb{\text{f2}[\text{x$\_$}]\text{:=}\text{f1}[x]-1-x}\)

\noindent\(\pmb{\text{f3}[\text{x$\_$}]\text{:=}\text{f2}[x]-1-x-\frac{x^2}{2}}\)

\noindent\(\pmb{\text{e1}[\text{x$\_$}]\text{:=}\frac{\text{f2}[x]}{\text{f3}[x]}-1}\)

\noindent\(\pmb{\text{e2}[\text{x$\_$}]\text{:=}\frac{\text{f0}[x]}{\text{f2}[x]}-1}\)

\noindent\(\pmb{\text{e3}[\text{x$\_$}]\text{:=}\frac{\text{f0}[x]}{\text{f3}[x]}-1}\)

\noindent\(\pmb{\text{e4}[\text{x$\_$}]\text{:=}\frac{\text{e1}[x]}{1+\text{e1}[x]}}\)

\noindent\(\pmb{\text{e5}[\text{x$\_$}]\text{:=}\frac{(1+\text{e2}[x])(1+\text{e3}[x])x^3}{8\text{f0}[x]}}\)

\noindent\(\pmb{\text{e6}[\text{x$\_$}]\text{:=}\frac{x^2(1+\text{e2}[x])^2}{\text{f0}[x]}}\)

\noindent\(\pmb{\text{e7}[\text{x$\_$}]\text{:=}\frac{\text{e4}[x]+\text{e5}[x]+\text{e6}[x]}{1-\text{e5}[x]}}\)

\noindent\(\pmb{\text{e8}[\text{x$\_$}]\text{:=}\frac{\text{e1}[x]+\text{e5}[x]}{1-\text{e5}[x]}}\)

\noindent\(\pmb{\text{e9}[\text{x$\_$}]\text{:=}x \text{e4}[x]}\)

\noindent\(\pmb{\text{e10}[\text{x$\_$}]\text{:=}\frac{2x \text{e4}[x]}{1-\text{e4}[x]}}\)

\noindent\(\pmb{\text{e11}[\text{x$\_$}]\text{:=}\frac{x(\text{e2}[x]+\text{e5}[x])+\text{e5}[x]}{(1-\text{e4}[x])(1-\text{e5}[x])}}\)

\noindent\(\pmb{d[\text{x$\_$}]\text{:=}\text{Max}[\text{e1}[x],\text{e2}[x],\text{e3}[x],\text{e4}[x],\text{e5}[x],\text{e6}[x],\text{e7}[x],\text{e8}[x]]}\)

\noindent\(\pmb{N[d[16]]}\)

\noindent\(0.000102752\)

\noindent\(\pmb{T=1-\frac{1}{2^{1/2}\text{Exp}[\text{Pi}/4]}}\)

\noindent\(1-\frac{e^{-\pi /4}}{\sqrt{2}}\)

\noindent\(\pmb{B=-.01+2(1-T)}\)

\noindent\(0.634794\)

\noindent\(\pmb{\text{A0}=B(\text{Exp}[2T/B]-1)}\)

\noindent\(4.73302\)

\noindent\(\pmb{\text{A1}=N\left[\frac{2(1+T)\text{Exp}[-\text{ArcTan}[T]]}{\left(1+T^2\right)^{1/2}}\right]}\)

\noindent\(1.5312\)

\noindent\(\pmb{\text{A2}=N\left[\frac{(1+T)\text{Exp}[-\text{ArcTan}[T]]}{\left(1+T^2\right)^{1/2}}\text{Integrate}\left[\frac{2\text{Exp}[\text{ArcTan}[x]]}{(1+x)\left(1+x^2\right)^{1/2}},\{x,0,1\}\right]\right]}\)

\noindent\(1.41846\)

\noindent\(\pmb{\text{A3}[\text{c$\_$},\text{x$\_$}]\text{:=}\frac{10 \text{A0} c }{B^2}+\frac{8 \text{A0}^2 c d[x]}{B^3}+\frac{(c+1)d[x]}{B(1-d[x])^2}}\)

\noindent\(\pmb{\text{A4}[\text{c$\_$},\text{x$\_$}]\text{:=}\text{A0} \text{A3}[c,x]/2}\)

\noindent\(\pmb{\text{A5}[\text{c$\_$},\text{x$\_$}]\text{:=}\frac{2c}{B}+\frac{2 \text{A4}[c,x]}{B}+\frac{4 c \text{A0}}{B^2}}\)

\noindent\(\pmb{\text{C1}[\text{c$\_$},\text{x$\_$}]\text{:=}\text{A1} c-\text{A2}-\text{A5}[c,x] d[x]}\)

\noindent\(\pmb{N[\text{C1}[15,16]]}\)

\noindent\(20.1217\)

\end{document}